\documentclass[11pt]{article}
\usepackage{bbm}
\usepackage[english]{babel}
\usepackage[latin1]{inputenc}
\usepackage{booktabs}
\usepackage{graphicx}
\usepackage{graphics}
\usepackage{caption}
\usepackage{subcaption}
\usepackage{amsmath}
\usepackage{amssymb}
\usepackage{amsthm}
\usepackage{hyperref}
\usepackage{array}
\usepackage{color}
 \usepackage[]{tabularx}

\usepackage[percent]{overpic}
\usepackage{float}

\theoremstyle{plain}% default
\newtheorem{theo}{Theorem}[section]
\newtheorem{lem}[theo]{Lemma}

\newtheorem{remark}[theo]{Remark}

\theoremstyle{definition}
\newtheorem{definition}{Definition}[section]

\newtheorem{example}{Example}[section]
\theoremstyle{remark}

\newcommand{\E}{\mathbb{E}}
\newcommand{\p}{\mathbb{P}}
\newcommand{\R}{\mathbb{R}}

\newcommand{\LL}{\mathcal{L}}

\newcommand{\PP}{\mathcal{P}}

\def\argmin{\textit{argmin}}
\title{Obtaining fairness using optimal transport theory}

\author{Eustasio del Barrio$^{a}$, Fabrice Gamboa$^{b}$, Paula Gordaliza$^{a,b}$,  and Jean-Michel 
Loubes$^{b}$
\\
$^{a}$\textit{IMUVA, Universidad de Valladolid} and $^{b}$\textit{Institut de Math\'ematiques de Toulouse}}

\numberwithin{equation}{section}

\begin{document}

\maketitle

\section*{Abstract}

Statistical algorithms are usually helping in making decisions in many aspects of our lives. But, how do we know if these algorithms are biased and commit unfair discrimination of a particular group of people, typically a minority? \textit{Fairness} is generally studied in a probabilistic framework where it is assumed that there exists a protected variable, whose use as an input of the algorithm may imply discrimination. There are different definitions of Fairness in the literature. In this paper we focus on two of them which are called Disparate Impact (DI) and Balanced Error Rate (BER). Both are based on the outcome of the algorithm across the different groups determined by the protected variable. The relationship between these two notions is also studied. The goals of this paper are to detect when a binary classification rule lacks fairness and to try to fight against the potential discrimination attributable to it. This can be done by modifying either the classifiers or the data itself. Our work falls into the second category and modifies the input data using optimal transport theory. 
\vskip .2in
{\bf Keywords :} Fairness in Machine Learning, Optimal Transport, Wasserstein barycenter. \vskip .2in

\section{Introduction}
Along the last decade, Machine Learning methods have become more popular to build decision algorithms. Originally meant for recommendation algorithms over the Internet, they are now widely used in a large number of  very sensitive areas such as medicine, human ressources with hiring policies, banks and insurance (lending), police, and  justice with criminal sentencing, see for instance in \cite{berk2017fairness} \cite{pedreschi2012study} or \cite{2018arXiv180204422F} and references therein. The decisions made by what is now referred to as IA have a growing impact on human's life.  The whole machinery of these technics relies on the fact that a decision rule can be learnt by looking at a set of labeled examples called the learning sample and then this decision will be applied for the whole population which is assumed to follow the same underlying distribution.  So the decision is highly influenced by the choice of the learning set.\\
 \indent In some cases, this learning sample may present some bias or discrimination that could possibly be learnt by the algorithm  and then propagated  to the entire population by automatic decisions and, even worse, providing a mathematical legitimacy for this unfair treatment.  When giving algorithms the power to make automatic decisions, the danger may come that the reality may be shaped according to their prediction, thus reinforcing their beliefs in the model which is learnt. Classification algorithms are one particular focus of fairness concerns since  classifiers map individuals to outcomes. \vskip .1in
Hence, achieving fair treatment  is one of the growing fields of interest in Machine Learning. We refer for instance to \cite{ZVGG} or \cite{2018arXiv180204422F} for a recent survey on this topic. For this, several definitions of fairness have been considered. In this paper we focus on the notion of disparate impact for protected variables introduced in \cite{FFMSV}. Actually, some variables, such as sex, age or ethnic origin, are potentially sources of  unfair treatment since they enable to create  information that should not be processed out by the algorithm. Such variables are called in the literature protected variables. An algorithm is called fair with respect to these attributes when its outcome does not allow to make inference on the information they convey.  Of course the naive solution of ignoring these attributes when learning the classifier does not ensure this, since the protected variables  may be closely correlated with other features enabling a classifier to reconstruct them. \\
\indent Two solutions have been considered in the Machine Learning literature. The first one consists in changing the classifier in order to make it not correlated to the protected attribute. We refer for instance to \cite{ZVGG}, \cite{2017arXiv170700044B} or \cite{2018arXiv180208626D} and references therein. Yet changing the way a model is built or explaining how the classifier is chosen may be seen too intrusive for many companies or some may not be able to change the way they build the model. Hence a second solution consists in changing the input data so that predictability of the protected attribute is impossible. The data will be blurred in order to obtain a fair treatment of the protected class. This point of view has been proposed in \cite{FFMSV}, \cite{johndrow2017algorithm}  or \cite{DBLP:journals/corr/abs-1712-07924} for instance. \vskip .1in
In the following, we first provide a statistical analysis of the Disparate Impact definition and recast some of the ideas developed in \cite{FFMSV} to stress the links between fairness, predictability and the distance between the distributions of the variables given the protected attribute. Then we provide some theoretical justifications of the methodology proposed by previous authors for one dimensional data to blur the data using the barycenter of the conditional distribution with respect to the Wasserstein distance. These methods are called  either  full or  partial repair. We extend this reparation procedure to the case of multidimensional data and provide a feasible algorithm to achieve this fairness reparation using the notion of  Wasserstein barycenter. Finally, we propose another methodology called {\it Random Repair} to transform the data in order to achieve a tradeoff between a minimal loss of information with respect to the classification task and still a certain level of fairness for classification procedures that could be used with this transformed data.  Applications to real data enable to study the efficiency of previous procedures.
\vskip .1in
The paper falls into the following parts. Section~\ref{s:recast} presents the relationships between the notions of Disparate Impact,  the predictability of a protected attribute and distance between the distributions conditionally to this attribute. Section~\ref{s:ot} is devoted to a probabilistic framework to transform the data to obtain fair classifiers. The following section, Section~\ref{s:repair}, provides some insight to understand the use of the Wasserstein's barycenter and its limitation. Applications to a real data set are shown in Section~\ref{s:num}, while the proofs are postponed to the Appendix. 
%Conclusions are given in Section~\ref{s:conclu}.

%Statistical algorithms are usually helping in making decisions in many aspects of our lives. But, how do we know if these algorithms are biased, involve ilegal discrimination or are unfair? \textit{Fairness} is generally studied in the framework where it is assumed that there exists a protected variable, whose use as an input of the algorithm may imply discrimination. In this paper, a review of the main definitions of fairness in Machine Learning is done. Then, we focus on two of them, which are based either in the outcome of the algorithm or in the error committed by it, across the different groups determined by the protected variable. When a procedure lacks of one of this two kinds of fairness, we will say that it has Disparate Impact (DI) and Disparate Mistreatment (DM), respectively. In these situations, our goal is both to detect and fight against discrimination. This can be done by modifying either the classifiers or the data itself. Our work falls into the second category and modifies the input data using optimal transport theory.
\section{Fairness using Disparate Impact assessment} \label{s:recast}
Consider the probability space $\left(\Omega, \mathcal{B}, \p\right)$, with $\mathcal{B}$ the Borel $\sigma-$algebra of subsets of $\R^d$ and $d \geqslant 1$. In this paper, we tackle the problem of forecasting a binary variable $Y: \Omega \rightarrow \left\lbrace0,1\right\rbrace$, using observed covariates  $X: \Omega \rightarrow \R^d, \ d\geqslant 1.$  We assume moreover that the population can be divided into two categories that represent a bias,  modeled by a  variable $S: \Omega \rightarrow \{0,1\}$. This variable is called the protected attribute, which takes the values $S=0$ for the \textquotedblleft minority\textquotedblright class and supposed to be the unfavored class; and $S=1$ for the \textquotedblleft default\textquotedblright, and usually favored class. We also introduce also a notion of positive prediction in the sense that $Y=1$ represents a success while $Y=0$ is a failure.  \\
Hence the classification problem aims at predicting a success from the variables $X$, using a family of binary classifiers $g\in \mathcal{G}: \R^d \rightarrow \{0,1\}$. For every $ g\in\mathcal{G}$, the outcome of the classification will be the prediction $\hat{Y}=g(X)$. We refer for instance to \cite{bousquet2004introduction} for a complete description of classification problems in statistical learning.  \vskip .1in
In this framework, discrimination or unfairness of the classification procedures, appears as soon as the prediction and the protected attribute are too closely related, in the sense that statistical inference on $Y$ may lead to learn the distribution of the protected attribute $S$.  This issue has received lots of interest among the last years and several ways to quantify this {\it discrimination bias} have been given. We highlight two of them, whose interest depends on the particular problem. More precisely, we can deal with two situations, depending whether the true distribution of the label $Y$ is available. If it is known,  the definition introduced in \cite{berk2017fairness}, defines that a classifier $g: \R^d \rightarrow \{0,1\}$ achieves \textit{Overall Accuracy Equality}, with respect to the joint distribution of $(X,S,Y)$, if
	\begin{equation} \label{OAE}
	\p(g(X)=Y \mid S=0) = \p(g(X)=Y \mid S=1).
	\end{equation}
	This entails that the probability of a correct classification is the same across groups and, hence, the classification error is independent of the group. This idea can be also found in the \cite{ZVGG} as the condition of $g$ having \textit{Disparate Mistreatment}, which happens when the probability of error is different for each group as in \eqref{OAE}. \\
	
\indent Nevertheless, in many problems, the true $Y$ is not available (this data may be very sensitive and the owner of the data may not want to make it available), or the classification methodology can not be changed, so the study of fairness must be based on the outcome $\hat{Y}$. In this situation, following \cite{FFMSV}  or \cite{berk2017fairness},  a classifier $g: \R^d \rightarrow \{0,1\}$ is said to achieve \textit{Statistical Parity}, with respect to the joint distribution of $(X,S)$, if
	\begin{equation} \label{SP} 
	\p(g(X)=1 \mid S=0) = \p(g(X)=1 \mid S=1).
	\end{equation}
	%we use a definition of \textit{Statistical Parity} \citep{Berk} in terms of the prediction $\hat{Y}=g(X)$.
	%\begin{definition}
%	A classifier $g: \R^d \rightarrow \{0,1\}$ achieves \textit{Statistical Parity}, with respect to the joint distribution of $(X,S)$, if
%		\begin{equation*}
%		\p(g(X)=1 \mid S=0) = \p(g(X)=1 \mid S=1).
%		\end{equation*}
%	\end{definition}
	This means that the probability of a successful outcome is the same across the groups. For instance, if we consider that the protected variable represents gender, the value $S=0$ would be assigned to \textquotedblleft female\textquotedblright and $S=1$ to \textquotedblleft male\textquotedblright, we would say that the algorithm used by a company achieves \textit{Statistical Parity} if a man and a woman have the same probability of success (for instance being hired or promoted). \\
	We will use the following notations 
	$$a(g):=\p(g(X)=1 \mid S=0), \quad b(g):=\p(g(X)=1 \mid S=1).$$ In the rest of the paper, we consider classifiers $g$ such that $a(g)>0$ and $b(g)>0$, which means that the classifier is not totally fair or unfair in the sense that the classifier does not predict the same outcome for a whole population according to the protected attribute. 
	
	\indent The independence described in \eqref{SP} is difficult to achieve and may not exist in real data. Therefore, to assess this kind of fairness, an index called \textit{Disparate Impact of the classifier $g$ with respect to $(X,S)$}, has been introduced in~\cite{FFMSV} as
	\begin{equation} \label{DI}
DI(g,X,S)=\frac{\p(g(X)=1 \mid S=0)}{\p(g(X)=1 \mid S=1)}.
\end{equation}
The ideal scenario where $g$ achieves \textit{Statistical Parity} is equivalent to $DI(g,X,S)=1$. Statistical Parity is often unrealistic, so we will relax it into  achieving a certain level of fairness as described in the following definition.

\begin{definition}
	The classifier $g: \R^d \rightarrow \{0,1\}$ has Disparate Impact at level $\tau \in \left(0,1\right]$, with respect to $(X,S)$,
 if $DI(g,X,S) \leq \tau$.
\end{definition}
Note the Disparate Impact of a classifier measures its level of fairness: the smaller the value of $\tau$, the less fair it is. The classification rules considered in this framework are such that $b(g) \geqslant a(g)>0,$ because we are assuming that the default class $S=1$ is more likely to have a successful outcome. Thus, in the definition, the level of fairness $\tau$ takes values $0 < \tau \leq 1$. We point out that the value $\tau_0=0.8=4/5$, which is also known in the literature as the \textit{$80 \%$ rule} has been cited as a legal score to decide whether the discrimination of the algorithm is acceptable or not (see for instance \cite{FFMSV}). This rule can be explained as \textquotedblleft for every 5 individuals with successful outcome in the majority class, 4 in the minority class will have a successful outcome too\textquotedblright.\\\indent  In what follows, to promote fairness, it will be useful to state the definition in the reverse sense. A classifier does not have Disparate Impact at level $\tau$, with respect to $(X,S)$, if $DI(g,X,S)>\tau$. \vskip .1in
Finally, another definition has been proposed in the statistical literature on fair learning. Given a classifier $g \in \mathcal{G}$, its Balanced Error Rate (BER) with respect to the joint distribution of the random vector $(X,S)$ is defined as the average class-conditional error
\begin{equation}\label{def:BER} 
BER(g,X,S)=\frac{\p\left(g(X)=0 \mid S=1\right)+\p\left(g(X)=1 \mid S=0\right)}{2}.
\end{equation}
Notice that $BER(g,X,S)$ is the general misclassification error of $g \in \mathcal{G}$ in the particular case when we have $\p(S=0)=\p(S=1)=1/2$, which consists in the ideal situation when both protected classes have the same probability of occurence. This quantity enables to define the notion of $\varepsilon-$predictability of the protected attribute.  $S$ is said to be $\varepsilon-$predictable from $X$ if there exists a classifier $g \in \mathcal{G}$ such that
	$$
	BER(g,X,S) \leq \varepsilon.
$$ Equivalently, $S$ is said not to be $\varepsilon-$predictable from $X$ if $BER(g,X,S) > \varepsilon$, for all classifiers  $ g$ chosen in the class $  \mathcal{G}$. Thus, if the minimum of this quantity is achieved by a classifier $g^*$,
$$ \min_{g \in \mathcal{G}}BER(g,X,S)= BER(g^*,X,S)=\varepsilon^*$$
then it is clear that $S$ is not $\varepsilon-$predictable from $X$ for all $ \varepsilon \leq \varepsilon^*$.\vskip .1in
In the following, we recast previous notions of fairness and provide a probabilistic framework to highlight the relationships between the distribution of the observations and the fairness of the classification problem. \vskip .1in
The following theorem generalizes the result in \cite{FFMSV}, showing the relationship between predictability and Disparate impact.

\begin{theo} \label{theo:DIlinkBER}
	Given random variables $X \in \R^d, \ S\in\{0,1\}$, the classifier $g \in \mathcal{G}$ has Disparate Impact at level $\tau \in \left[0,1\right]$, with respect to $(X,S)$, if and only if $BER(g,X,S) \leq \frac{1}{2}-\frac{a(g)}{2}(\frac{1}{\tau}-1)$.
\end{theo}
% In the paper where this theorem is presented, the level of DI is taken equal to $\tau_0=0.8$. In this case,
%		\begin{center}
%		a classifier $g:X\rightarrow Y$ has DI at level $0.8$ if, and only if, $S$ is $\left(\frac{1}{2}-\frac{\alpha(g)}{8}\right)-$predictable from $X$.
%		\end{center}
 
% Now, we will see that if we are able to predict the protected variable $S$ from $X$, then we will have Disparate Impact at some level $\tau \in \left[0,1\right]$. To show this claim, we shall introduce the notion of \textit{predictability} of $S$ from $X$, which is based on a particular measure of the loss in the prediction of $S$ from $X$, which is called \textit{Balanced Error Rate} \citep{FFMSV}.

 The following theorem establishes the relationship between $\varepsilon^*$ the minimum Balance Error Rate and  distance in Total Variation between the two conditional distributions $\mathcal{L}\left(X|S=0\right)$ and $\mathcal{L}\left(X|S=1\right)$.

\begin{theo}\label{theo:minBER}
	Given the variables $X: \Omega \rightarrow \R^d, \ d \geqslant 1,$ and $S: \Omega \rightarrow \{0,1\}$,
	\begin{equation*}
	\min_{g \in \mathcal{G}}BER(g,X,S)=\frac{1}{2}\left(1-d_{TV}\left(\mathcal{L}\left(X|S=0\right), \mathcal{L}\left(X|S=1\right)\right)\right).
	\end{equation*}
%	where $g: \R^d \rightarrow \{0,1\}$ varies in the set of binary classifiers $\mathcal{G}$.
\end{theo}

This result shows that fairness expressed through the notion of Disparate Impact depends highly on the conditional distributions of the variables X conditionally to the protected attribute, $\mathcal{L}\left(X|S=0\right)$ and $ \mathcal{L}\left(X|S=1\right)$. \vskip .1in

Actually, Theorem~\ref{theo:minBER} implies that  $S$ is not $\varepsilon-$predictable from $X$ if, and only if,
\begin{equation}\label{equ:equivSnotPREDandDTV}
d_{TV}\left(\mathcal{L}\left(X|S=0\right), \mathcal{L}\left(X|S=1\right)\right) < 1-2\varepsilon.
\end{equation}
and, as a consequence of Theorem~\ref{theo:DIlinkBER}, for all $g \in \mathcal{G},$
\begin{equation*}
DI(g,X,S)>\displaystyle\frac{1}{1+\displaystyle\frac{d_{TV}\left(\mathcal{L}\left(X|S=0\right), \mathcal{L}\left(X|S=1\right)\right)}{a(g)}}.
\end{equation*}
Hence, the smaller the Total Variation distance, the greater is the value $\varepsilon$ that we could find satisfying Equation (\ref{equ:equivSnotPREDandDTV}) and thus, the less predictable using the variables $X$ will be $S$ . The best case happens when $d_{TV}\left(\mathcal{L}\left(X|S=0\right), \mathcal{L}\left(X|S=1\right)\right)=0$, which is equivalent to the equality of both conditional distributions $\mathcal{L}\left(X|S=0\right) = \mathcal{L}\left(X|S=1\right)$. In this situation, $X$ and $S$ are independent random variables, and we will have that $S$ is not $\varepsilon-$predictable from $X$, for all  $\varepsilon \leqslant \frac{1}{2}$, and $DI(g,X,S)=1$.  Note that clearly $\varepsilon = 1/2$ non predictability is the best that can be achieved. \vskip .1in

\section{Removing disparate impact using Optimal Transport} \label{s:ot}
\subsection{A probabilistic model for data repair} 
Some classification procedures exhibit a discrimination bias quantified through a potential Disparate Impact in the classification outcome $\hat{Y}=g(X)$, with respect to the joint distribution of $(X,S)$. To get rid of the possible discrimination associated to a classifier $g$, two main strategies can be used, either modifying the classifiers or modifying the input data. In this work, we are facing the problem where we have no access to the values $Y$ of the learning sample, hence we focus on the methodologies that intend to modify the data in order to achieve fairness. \\
The main idea is to change the data in order to break their relationship with  the protected attribute. This transformation is called repairing the data. For this, \cite{FFMSV}, \cite{johndrow2017algorithm}  or \cite{DBLP:journals/corr/abs-1712-07924} propose to map the conditional distributions to a common distribution in order to achieve statistical parity as described in \eqref{SP}. The choice of the common distribution in one dimension is described as the distribution obtained by taking the mean of the quantile functions.  A total repair of the data amounts to modifying the input variables $X$ building a repaired version, denoted by $\tilde{X}$, such that any classifier $g$ will have Disparate Impact at level $\tau=1$, with respect to $(\tilde{X},S)$. This means that every classifier $g$ used to predict the target class $Y$ from the new variable $\tilde{X}$ will achieve \textit{Statistical Parity} with respect to $(\tilde{X},S)$. As a counterpart, it is clear that the choice of the distribution to whom the original  variables are mapped should convey as much as information possible on the original variable, otherwise it would hamper the accuracy  of the new classification. This constraint led some authors to recommend the use of the so-called Wasserstein barycenter. \vskip .1in We now present some statistical justifications for this choice and provide some comments on the way to repair the data to obtain fair enough classification rules without modifying too much the original data set.

Achieving Statistical Parity amounts to modifying the original data into a new random variable $\tilde{X}$  such that the conditional distribution with respect to the protected attribute $S$ is the same for all groups, namely 
\begin{equation}\label{eq:fairnessrequirement}
\mathcal{L}\left(\tilde{X} \mid S=0\right) = \mathcal{L}\left(\tilde{X}\mid S=1\right).
\end{equation}
In this case, any classifier $g$ built with such information will be such that  
$$\mathcal{L}\left(g(\tilde{X})\mid S=0\right) = \mathcal{L}\left(g(\tilde{X})\mid S=1\right),$$ which implies that $DI(g, \tilde{X},S)=1$ and so this transformation promotes full fairness of the classification rule. \vskip .1in
To achieve this transformation, the solution detailed in many papers is to map both conditional distributions $\mu_0 := \mathcal{L}(X|S=0)$ and $\mu_1 := \mathcal{L}(X|S=1)$ onto a common distribution $\nu$. Actually, the distribution of the original variables $X$ is transformed using a map $T_S$ which depends on the value of the protected attribute $S$

%Formally, we want to find the function
\begin{center}
	$
	\begin{array}{cccc}
	T_S: & \R^d & \longrightarrow  & \R^d  \\
	& X  & \longmapsto & \tilde{X}=T_S(X),
	\end{array}
	$
\end{center}
and such that 
\begin{equation}\label{equ:conditionfairtransf}
\mathcal{L}\left(T_0(X)\mid S=0\right) = \mathcal{L}\left(T_1(X)\mid S=1\right).
\end{equation}
Note that the function $T_S$ is random because of its dependence on the binary random variable $S$. \\

\begin{figure}[h]
	\centering
	\begin{overpic}[scale=.5,unit=1mm]{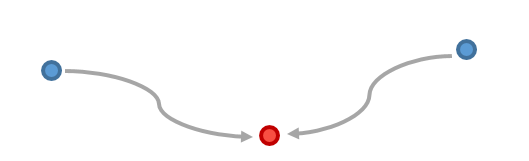}
		\put (6,22){$\mu_0$}
		\put (87,26){$\mu_1$}
		\put (50,10){$\nu$}
	\end{overpic}
	\caption{General repairing scheme}
	\label{fig:push}
\end{figure}

In this framework, the problem of achieving Statistical Parity is the same as the problem of finding a (random) function $T_S$ such that (\ref{equ:conditionfairtransf}) holds. As it is represented in Figure \ref{fig:push}, if we denote by $\mu_S \sim X \mid S$, our goal is to map these two distributions  to a common law $\nu=\mu_S \circ T^{-1}_S$.
Consequently, two different problems arise
\begin{itemize}
\item First of all,  the choice of the distribution $\nu$ should be as similar as possible to both distributions $\mu_0$ and $\mu_1$ at the same time, in order to reduce the amount of information lost with this transformation and thus still enabling the  prediction task  using the modified variable $\tilde{X} \sim \nu$  instead of the original $X$.
\item On the other hand, once we have selected the distribution $\nu$, we have to find the optimal way of transporting $\mu_1$ and $\mu_0$  to this new distribution $\nu$. 
\end{itemize}
From Section~\ref{s:recast}, the natural distance related to fairness between the two conditional distributions  is the total variation distance and that should be used. However, this distance is computationally difficult to handle, hence previous works promote the use of Wasserstein distance which appears as a natural distance to move distributions. \vskip .1in

For this, we recall some results on optimal transport theory and Wasserstein metric between probability measures, which provides an appropriate tool for comparing probability distributions. In this framework, the map $T_S$ will be a random transport plan between the distributions $\mathcal{L}(X\mid S)$ and $\mathcal{L}(\tilde{X})$. Moreover, we will first recall the definition of Wasserstein barycenters which are often chosen in the statistical literature as   new distribution $\nu$. 

\subsection{Wasserstein distance and Wasserstein barycenters}
Consider the space $\PP_2\equiv\PP_2(\R^d)$ of Borel probabilities on $\R^d$ with finite second moment. The related set $\PP_{2,ac}\equiv\PP_{2,ac}(\R^d)$ will denote the subset of $\PP_2(\R^d)$ containing the probabilities that are absolutely continuous with respect to Lebesgue measure. Given $\mu, \nu \in \PP_2(\R^d)$,  we denote by $\Pi(\mu, \nu)$ the set of all
probability measures $\pi$ over the product set $\R^d \times\R^d $
with first (resp. second) marginal $\mu$ (resp. $\nu$).

The transportation cost with quadratic cost function, or quadratic transportation cost,
between these two measures
$\mu$ and $\nu$   is defined as
\begin{equation*}
\label{equ:distWass}
 \mathcal{T}_2(\mu, \nu) = \inf_{ \pi \in \Pi(\mu, \nu)} \int \left\Vert x - y\right\Vert ^2 d \pi(x,y).
\end{equation*}
The quadratic transportation cost allows to endow the set $\PP_2$  with a metric by defining  the $2$-Wasserstein distance between $\mu$ and $\nu$ as
%\begin{equation*}
$ W_2(\mu, \nu)  = \mathcal{T}_2(\mu, \nu) ^{1/2}.$
%\end{equation*}
More details on Wasserstein distances and their links with optimal transport problems can be found in  
\cite{rachev} or \cite{villani2008optimal} for instance.

In the one-dimensional case $W_{2}(\mu, \nu)$ is simply the $L_2-$distance between the quantile functions of $\mu$ and $\nu$, enabling  its direct computation
\begin{equation*}
W_{2}^2(\mu, \nu)= \int_{0}^1\left|F^{-1}(t) - G^{-1}(t)\right|^2dt, \ F \sim \mu, \ G \sim \nu.
\end{equation*}

 A distribution $\pi$ with marginals $\mu$ and $\nu$ which minimizes  \eqref{equ:distWass} is called an optimal coupling of $\mu$ and $\nu$. Moreover, if $\mu$ vanishes on sets of dimension $d-1$, in particular if $\mu \in \PP_{2,ac}(\R^d)$, then there exists an optimal transport map, $T$, transporting (pushing forward) $\mu$ to $\nu$. The following Theorem is a convenient version that can be found in \cite[Theorem 2.12]{Villani2003} .

\begin{theo}\label{theo:optimalcoupling}
	Let $\mu, \nu \in \PP_2(\R^d)$ and let $\pi=\LL(X,Y)$ be the joint distribution of a pair $(X,Y)$ of $R^d-$valued random vectors with probability laws $\LL(X)=\mu$ and $\LL(Y)=\nu.$
	\begin{enumerate}
		\item[(i)]The probability distribution $\pi$ is an optimal coupling of $\mu$ and $\nu$ if, and only if, there exists a convex lower semi-continuous function $\varphi \in \partial\varphi(X)$ a.s. such that $\pi$ is concentrated on $\partial \varphi,$ the subgradient of $\varphi$.
		\item[(ii)]If we assume that $\mu$ does not give mass to sets of dimension at most $d-1$, then there is a unique optimal coupling $\pi$ of $\mu$ and $\nu$, that can be characterized as the unique solution to the Monge transportation problem - an optimal transport map - $T$, i.e.: $\pi= \mu \circ (Id,T)^{-1}$ (or $Y=T(X)$ a.s.), and
		\begin{align*}
		W_2^2(\mu,\nu) & =\int\left\|x-T(x)\right\|^2d\mu(x)\\
		& =\inf\left\{\int\left\|x-S(x)\right\|^2d\mu(x), \ where \ S \ satisfies \ \nu=\mu\circ S^{-1}\right\}
		\end{align*}
		Such a map is characterized as $T=\nabla\varphi \ \mu-a.s.,$ the $\mu-a.s.$ unique function that maps $\mu$ to $\nu$ and that is the gradient of a lower semi-continuous funtion $\varphi$. In the following we will use the notation $$ \nu=T_\sharp \mu = \mu \circ T^{-1}.$$
			\end{enumerate}
\end{theo}

We point out that the existence of the o.t map  is commonly referred to as Brenier's theorem and originated from Y. Brenier's work in the analysis and mechanics literature. However, it is worthwile pointing out that a similar statement was established earlier independently in a
probabilistic framework by J.A. Cuesta-Albertos and C. Matr\'{a}n \cite{cuesta1989notes} : they show existence of an optimal transport map for quadratic cost over Euclidean and Hilbert spaces, and prove monotonicity of the optimal map in some sense (Zarantarello monotonicity). \vskip .1in

When dealing with a collection of distributions $\mu_1,\dots,\mu_J$,  we can define a notion of variation of these distributions. For any $\nu \in \PP_2(\R^d)$, set 
$$ {V_2}(\nu)=\sum_{j=1}^J  \omega_j W_2^2(\nu,\mu_j)$$
where  $\omega_1, \ldots, \omega_J$ are positive real numbers such that $\sum_{j=1}^J\omega_j=1$.
 This quantity provides a global measure of separation between the probabilities $\mu_j, \: j=1,\ldots,J,$ with respect to fixed weights and has been received recently.  Finding the distribution minimizing the variance of the distributions has been tackled when defining the notion of barycenter of distributions with respect to Wasserstein's distance in the seminal work of~\cite{agueh2011barycenters}.  More precisely, given $p\geq 1$, they provide conditions to ensure existence and uniqueness of the barycenter of the probability measures $(\mu_j)_{1\leq j \leq J}$  with weights $(\omega_j)_{1\leq j \leq J}$, i.e. a minimizer of the following criterion
\begin{equation} \label{carlier}
\nu \mapsto \sum_{j=1}^J \omega_j W_2^2(\nu,\mu_j).
\end{equation}
 Such a minimizer, $\mu_B$, is called a barycenter or Fréchet mean of $\mu_1, \ldots, \mu_J$,  with respect to the weights $\omega_1,\ldots,\omega_J$. Empirical versions of the barycenter and their properties are analyzed in \cite{boissard2015distribution} or \cite{le2017existence}. Similar ideas have also been developed in \cite{cuturi2014fast} or \cite{bigot2012characterization}. Hence the Wasserstein barycenter  distribution appears to be a meaningful feature to represent the mean variations of a set of distributions. \vskip .1in We point out that its computation is a difficult issue for the general case. Yet, in this work, we only consider barycenter between two probabilities $\mu_0, \mu_1 \in \PP_2(\R^d)$.  For the one dimensional case, the solution proposed in \cite{FFMSV} to repair the data is to map these distributions to a distribution whose quantile function is defined by taking the mean of the quantile functions of $\mu_0$ and $\mu_1$. This corresponds actually to the minimizer of \eqref{carlier} for distributions on the real line denoted by $\PP_{2,ac}(\R)$. In the following, we present in Section~\ref{s:num} how to compute a barycenter between two distributions in higher dimensions and propose in Section~\ref{s:repair} a justification for using the Wasserstein barycenter to repair the data.

\section{Full and Partial Repair with Wasserstein Barycenter} \label{s:repair}

In our particular problem, where we have $J=2$, the two conditional distributions of the random variable $X$ by the protected attribute $S$ are going to be transformed into the distribution of the Wasserstein barycenter $\mu_B$ between $\mu_0$ and $\mu_1$, with weights $\pi_0$ and $\pi_1$, defined as 
\begin{equation*}
\mu_B \in \argmin_{\nu \in \PP_2}V^2_2(\mu_0,\mu_1;\pi_0,\pi_1)=\argmin_{\nu \in \PP_2}  \left\{\pi_0 W_{2}^2(\mu_0, \nu) + \pi_1W_{2}^2(\mu_1, \nu)\right\}
.
\end{equation*}
Let $\tilde{X}$ be the transformed variable with distribution $\mu_B$. For each $S=s$, the deformation will be performed through the optimal transport map $T_s : \R^d \rightarrow \R^d$ pushing each $\mu_s$ towards the weighted barycenter $\mu_B$, whose existence is guaranteed as soon as $\mu_s$ are absolutely continuous with respect to Lebesgue measure using Theorem~\ref{theo:optimalcoupling}, which satisfies

\begin{equation}\label{equ:Mongeproblem}
 \E\left(\left\|X-T_s(X)\right\|^2 \mid S=s\right)=W_2^2(\mu_s, \mu_B).
\end{equation}
\begin{remark}\label{remark:barycenter2measures}
 Note first that in the particular setting of two distributions, the computation of the barycenter of the two measures is equivalent to the computation of the optimal transport map between them. More precisely, if we assume that $\mu_0 \in \PP_{2,ac}(\R^d)$ and denote by $T: \R^d \rightarrow \R^d$ the optimal transport map between $\mu_0$ and $\mu_1$, that is $\mu_1={\mu_0}_\sharp T$, then we can write
\begin{equation*}
\mu_\lambda=\mu_0\sharp\left((1-\lambda)Id + \lambda T\right),
\end{equation*}
where the map $(1-\lambda)Id + \lambda T$ is an optimal transport plan, for all $\lambda \in \left[0,1\right]$. We have that the measure $\mu_\lambda$ is the weighted barycenter between $\mu_0$ and $\mu_1$, with weights $1-\lambda$ and $\lambda$, respectively. So, the complexity of computing $\mu_B=\mu_0\sharp\left(\pi_0Id + \pi_1T\right)$ is the same as the complexity of computing $T$. 
\end{remark}

\begin{remark} Note also that for distributions on the real line,  we can write the explicit expression of the barycenter $\mu_B$ based on  the exact solution to the optimization problem (\ref{equ:Mongeproblem}). Given $S=s$ and $X \in \R$, let $F_s : \R \rightarrow [0,1]$ denote the cumulative distribution function of $X$ given that $S=s$, and $F^{-1}_{s} : [0,1] \rightarrow \R$ its quantile associated function. The weighted Wasserstein barycenter $\mu_B$ of the two distributions $\mu_0$ and $\mu_1$ is the unique minimizer of the functional (\ref{carlier}) and its quantile function can be computed as
%	\begin{equation*}
%	F_B^{-1}(t)=\frac{1}{2}\left(F^{-1}_{0}(t) + F^{-1}_{1}(t)\right), \ t \in [0,1].
%	\end{equation*}
	\begin{equation*}
	F_B^{-1}(t)=\left(\lambda F^{-1}_{0}(t) + (1-\lambda)F^{-1}_{1}(t)\right), \ t \in [0,1].
	\end{equation*}
	Moreover, we note that
	\begin{equation*}
	F_s\left(X \mid S=s\right) \sim \mathcal{U}(0,1), \ s=0,1,
	\end{equation*}
	and the optimal transport map solution to (\ref{equ:Mongeproblem}) is $T_s=F_B^{-1} \circ F_s$.
\end{remark}

\begin{figure}[h]
	\centering
	\begin{overpic}[scale=.5,unit=1mm]{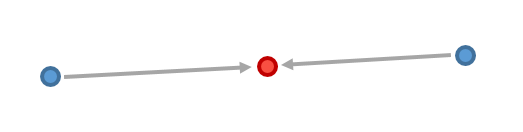}
		\put (5,15){$\mu_0$}
		\put (88,20){$\mu_1$}
		\put (25,15){$T_0$}
		\put (69,18){$T_1$}
		\put (48,18){$\mu_B$}
	\end{overpic}
	\caption{Repairing scheme towards the Wassersein barycenter}
	\label{fig:pushB}
\end{figure}
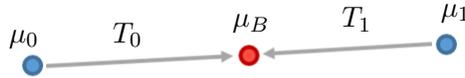
\subsection{Total repair}

To understand the use of  the Wasserstein barycenter distribution as the target distribution for $\mu_0$ and $\mu_1$, we quantify the amount of information lost by replacing the distribution of $X$ by the distribution of $\tilde{X}$ obtained by transporting these two distributions. Set the random transport plan $T_S: \R^d \longrightarrow  \R^d$, and the modified variable $\tilde{X}=T_S(X)$. We point out that choosing the distribution of $\tilde{X}$ amounts to choose the transportation plans $T_0$ and $T_1$. \vskip .1in

We are facing  learning problems in two different settings.
\begin{itemize}
\item On the one hand, the full information available is the input variables $X$ and also the protected variables $S$ which play an important role in the classification, since the classifier has a different behavior according to the class $S=0$ and $S=1$. Hence we let $S$ play a role in the decision process since it is associated to $Y$, and possibly giving rise to a different treatment for the two different groups. In this case, the classification risk when the full data $(X,S)$ is available can be computed as
 $R(g,X,S)$, the risk in the prediction of a classification rule $g $,  that depends on both variables $X$ and $S$, namely
\begin{equation*}
R(g,X,S):=\p(g(X,S)\neq Y).
\end{equation*}
\item On the other hand, in the repair data, only the modified version of the input data is at hand, $\tilde{X}.$ Hence learning a classifier amounts to minimizing 
\begin{equation*}
R(h,\tilde{X}):=\p(h(\tilde{X})\neq Y).
\end{equation*}
\end{itemize}
Studying the efficiency of the method requires providing a bound for the difference between the minimal risks obtained for the best classifier with input data $\tilde{X}=T_S(X)$, and for the best classifier with input data $(X,S),$ called $g_B$. These risks are respectively denoted $R_B(\tilde{X})$ and $R_B(X,S)=\inf_g R(g,X,S)=R(g_B,X,S),$ and then its difference is 
$$\mathcal{E}(\tilde{X}):=  R_B(\tilde{X}) - R_B(X,S).$$
Note first that, given $X=x$ and $S=s$, $\inf_g R(g,X,S)$ can be computed by mimicking the usual expression of the 2-class classification error as in \cite{bousquet2004introduction} for instance. We obtain
\begin{align*}
& \p(g(X,S)\neq Y \mid X=x, S=s)  \\
 = & \p(g(x,s)\neq 0, Y=0 \mid X=x, S=s) + \p(g(x,s)\neq 1, Y=1 \mid X=x, S=s)\\
= & \mathbbm{1}_{g(x,s) \neq 0} \p(Y=0 \mid X=x, S=s)+ \mathbbm{1}_{g(X) \neq 1}\p(Y=1 \mid X=x, S=s).
\end{align*}
Denoting the conditional expectation as
\begin{equation} \label{eq:cond}
\eta_s(x)=\p(Y=1 \mid X=x, S=s),
\end{equation}
we can write that
\begin{align*}
\p(g(X,S)\neq Y \mid X=x, S=s) & =(1-\mathbbm{1}_{g(x,s)=0})(1-\eta_s(x)) + \mathbbm{1}_{g(x,s) = 0} \eta_s(x)\\
& = \mathbbm{1}_{g(x,s) =0}(2\eta_s(x) -1)+ 1-\eta_s(x).
\end{align*}
Finally,  we get  \begin{equation} \label{eq:gen}
R(g,X,S)=\E \left[ \mathbbm{1}_{g(X,S)=0}(2\eta_S(X) -1)\right]+ \E \left[1-\eta_S(X) \right]. \end{equation}
 The minimum risk is thus obtained using  the Bayes' rule $g_B(x,s)= \mathbbm{1}_{\eta_s(x) > 1/2},$ leading to 
$$R_B(X,S):=inf_{g}R(g,X,S)=\E \left[ \mathbbm{1}_{\{2\eta_S(X) -1 <0\}}(2\eta_S(X) -1)\right]+ \E \left[1-\eta_S(X) \right].$$
\indent Similarly, the risk related to a classification rule $h(\tilde{X})$  is given by
\begin{equation}\label{equ:riskmodified}
R(h, \tilde{X})= R(h, T_S(X))=\E \left[ \mathbbm{1}_{h\circ T_S(X)=0}(2\eta_S(X) -1)\right]+ \E \left[1-\eta_S(X) \right].
\end{equation}
Hence, the amount of information lost due to modifying the data is controlled by the following theorem.

\begin{theo}\label{theo:bounddiffrisk}
	Consider  $X \in  \R^d$ and $S \in \{0,1\}$. Let $T_S:\R^d \rightarrow \R^d, \ d\geqslant 1$ be a random transformation of $X$ such that $\mathcal{L}(T_0(X) \mid S=0)=\mathcal{L}(T_1(X) \mid S=1)$, and consider the transformed version $\tilde{X}=T_S(X)$. For each $s \in \{0,1\}$, assume that the function $\eta_s(X)$ defined in \eqref{eq:cond} is Lipschitz with constant $K_s >0$
	Then, if $K=\max\{K_0,K_1\}$,
	\begin{equation}\label{equ:boundlossBayes}
	\mathcal{E}(\tilde{X}) \leq 2 \sqrt{2}K \left(\sum_{s=0,1}\pi_sW_2^2(\mu_s, {\mu_s}_\sharp T_s)\right)^\frac{1}{2}.
	\end{equation}
\end{theo}

The proof of this theorem which relies on the following lemma is postponed to the Appendix.
\begin{lem}\label{lem:bounddiffrisk}
Under Assumptions of Theorem~\ref{theo:bounddiffrisk}, the following bound holds
\begin{equation*}
R(g_B\circ T_S,X)-R(g_B,X,S)
\leq 2\E\left[\left|\eta_S(X)- \eta_S \circ T_S(X)\right|\right].
\end{equation*}
\end{lem}

%\begin{remark}
%	From (\ref{equ:riskmodified}), the risk in the classification with the modified variable when the classifier does not depend on the variable $S$ can be written as 
%	\begin{equation*}
%	R(g, \tilde{X})= R(g, T_S(X))=\E \left[ \mathbbm{1}_{g\circ T_S(X)=0}(2\eta_S(X) -1)\right]+ \E \left[1-\eta_S(X) \right],
%	\end{equation*}
%	and the lemma gives us
%	\begin{equation*}
%	R(g, T_S(X)) - R(g_B,X,S) \leq R(g_B\circ T_S,X,S)-R(g_B,X,S)
%	\leq 2\E\left[\left|\eta_S(X)- \eta_S \circ T_S(X)\right|\right].
%	\end{equation*}
%\end{remark}

%$R_B(\tilde{X}) \leq R(g_B, \tilde{X})$, and by Theorem \ref{theo:bounddiffrisk}
%	\begin{equation*}
%	R_B(\tilde{X}) - R_B(X,S)  \leq R(g_B, \tilde{X}) - R_B(X,S)  \leq 2 \sqrt{2}K \left(\sum_{s=0,1}\pi_sW_2^2(\mu_s, {\mu_s}_\sharp T_s)\right)^\frac{1}{2}.
%	\end{equation*}

Hence, Theorem \ref{theo:bounddiffrisk} provides some justification to the use of the Wasserstein barycenter as the distribution of the modified variable. Actually, minimizing the upper bound in~\eqref{equ:boundlossBayes} with respect to the function $T_S:\R^d \rightarrow \R^d, \ d\geqslant 1,$ leads to consider  the transport plan carrying the conditional distributions $\mu_s, \ s=0,1,$ towards their Wasserstein barycenter $\mu_B$ with weights $\pi_0, \pi_1$, that is, ${\mu_S}_\sharp T_S=\mu_B$. Hence, this provides some understanding on the choice of the Wasserstein barycenter advocated in the work \cite{FFMSV}. This leads to the following bound 
	\begin{align*}\label{equ:boundriskbarycenter}
\inf_{T_S}\left\{R(g_B\circ T_S,X)-R(g_B,X,S)\right\} & \leq 2 \sqrt{2}K \inf_{T_S} \left(\sum_{s=0,1}\pi_sW_2^2(\mu_s, {\mu_s}_\sharp T_s)\right)^\frac{1}{2}\\
& =\sqrt{2}K \left(\sum_{s=0,1}\pi_sW_2^2(\mu_s, \mu_B)\right)^\frac{1}{2}  \\
& = \sqrt{2}K \pi_0\pi_1 W_2^2(\mu_0, \mu_1)
% 2 \sqrt{2}KV_2\left(\mu_0, \mu_1; \pi_0, \pi_1\right).
	\end{align*}
Yet this bound is only an upper bound which only provides some guidelines on the choice of the choice of the distribution to which the conditional distribution have to be mapped.  Nevertheless, choosing the Wasserstein barycenter provides a simply a reasonable and, more important, feasible solution for fairness, to achieve statistical parity.

\subsection{Partial repair}\label{sec:partialrepair}

As pointed out in previous section, the Total Repair process  ensures full fairness but at the expense  of  the accuracy of the classification. A solution for this could be found in \cite{FFMSV},  called \textit{Geometric Repair}. The authors propose not to move all conditional distributions to the barycenter but only towards the barycenter on the Wasserstein's geodesic path  between $\mu_0$ and $\mu_1$. We analyze next this procedure and propose an alternative to this choice. \vskip .1in

Let $\lambda \in [0,1]$ be the parameter representing the amount of repair desired for $X$.  Let  $Z$ be a target variable with distribution $\mu$. Set $R_s=T_s^{-1}, \ s=0,1,$ where $T_s$ is the optimal transport map pushing each $\mu_s$ towards the target  $\mu$. In the literature, $\mu$ is chosen as the barycenter $\mu_B$ and the \textit{Partially Repaired} conditional distributions for $s\in \{0,1\}$ are defined as 
\begin{equation*}
	\mu_{s,\lambda}= \mathcal{L}(\lambda Z+ (1-\lambda)R_s(Z))=\mathcal{L}(\lambda T_s(X) + (1-\lambda)X \mid S=s).
\end{equation*}
\begin{figure}[H]
	\centering
	\begin{overpic}[scale=.5,unit=1mm]{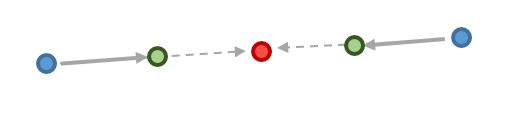}
		\put (6,11){$\mu_0$}
		\put (28,12){$\mu_{0,\lambda}$}
		\put (87,15){$\mu_1$}
		\put (67,13){$\mu_{1,\lambda}$}
		\put (48,13){$\mu_B$}
	\end{overpic}
	\caption{Original distributions (blue) and their partially repaired versions (green) towards the barycenter (red)}
	\label{fig:pushpartial}
\end{figure}
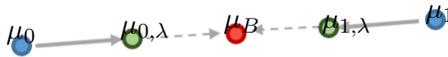
% Consider the geodesic functions from $\mu_S$ to $\mu_B$:
%\begin{equation}\label{equ:geodesic0tobarycenter}
%T_{0,\lambda}:=(1-\lambda\pi_1)Id+\lambda\pi_1 T,
%\end{equation}
%\begin{equation}\label{equ:geodesic1tobarycenter}
%T_{1,\lambda}:=(1-\lambda\pi_0)Id+\lambda\pi_0 T^{-1},
%\end{equation}
%where $T$ is the optimal transport map (o.t.m.) from $\mu_0$ towards $\mu_1$ (in the unidimensional setting $d=1$, $T=F_1^{-1}\circ F_0$), and take
%\begin{equation*}
%\tilde{X}_{S, \lambda}:=T_{S,\lambda}(X)
%\end{equation*}
%\begin{equation*}
%\mu_{S,\lambda}=T_{S,\lambda}\sharp\mu_S
%\end{equation*}
%the partially modified version of the original variable $X$ and its distribution, respectively.
%We note that for $\lambda=0$, the geodesic is the identity function $T_{S,0}=Id$, and for $\lambda_s=1$, it is the complete transport $T_{S,1}=T_S$. In this setting, as illustrated in Figure \ref{fig:pushpartial}, we will be creating a different distribution $\mu_{s,\lambda} = \mu_s \sharp T_{s,\lambda}, \ \forall \lambda \in [0,1]$.
This procedure is represented in Figure~\ref{fig:pushpartial}. Observe that $\lambda=1$ yields the fully repaired variable, and $\lambda=0$ leaves the conditional distributions unchanged. So the parameter $\lambda$ governs how close the distributions are to the barycenter.  Choosing the parameter $\lambda$ should be a trade-off between, on the one hand, accuracy of the classification error that leads to little change in the initial distribution, and, on the other hand, non predictability of the protected variable which implies that the two conditional distribution should be close with respect to the total variation distance.  \vskip .1in

Arguing among the lines of previous section to obtain an upper bound for the classification risk using the two distributions $\mu_{s,\lambda}, \ s\in \{0,1\},$ does not lead to a satisfying result. This comes from the fact that we move distributions according to Wasserstein distance, while fairness is measured using the total variation distance and they are of different nature. In fact, the distance in total variation between two probabilities $P$ and $Q$ can be computed as 
\begin{equation*}
d_{TV}(P,Q)=\min_{\pi \in \Pi(P,Q)}\pi(x\neq y),
\end{equation*}
see, e.g. \cite{massart2007concentration}.

So if $\lambda \in (0,1)$, this implies that
\begin{equation}\label{eq:bounddTVpartialrepair}
d_{TV}(\mu_{0,\lambda}, \mu_{1,\lambda}) \leq \p(\lambda Z + (1-\lambda)R_0(Z) \neq \lambda Z + (1-\lambda)R_1(Z)) = \p(R_0(Z) \neq R_1(Z)).
\end{equation}
Previous bound means that the amount of repair quantified by the parameter $\lambda$ does not affect the distance in Total Variation between the modified conditional distributions. Moreover, in some situations, \eqref{eq:bounddTVpartialrepair} turns out to be an equality. Consider, for instance, 
$$\mu_{0,0}=U(K,K+1)$$
$$\mu_{1,0}=U(-K-1,-K)$$
as the distributions of $X$ in each class. Then, the barycenter is $\mu_{0,1}=\mu_{1,1}=U(-1/2,1/2)$ and the partially repaired distributions are
$$\mu_{0,\lambda}=U(-\lambda/2+(1-\lambda)K,-\lambda/2+(1-\lambda)K+1)$$
$$\mu_{1,\lambda}=U(-\lambda/2-(1-\lambda)(K+1),-\lambda/2-(1-\lambda)(K+1)+1).$$

In this particular case, the distance in total variation can be easily computed as
\begin{equation*}
d_{TV}(\mu_{0,\lambda},\mu_{1,\lambda})=\min(1,(1-\lambda)(2K+1)).
\end{equation*}
As a consequence, $d_{TV}(\mu_{0,\lambda},\mu_{1,\lambda})=1,$ if $\lambda\leq
2K/(2K+1)$, which means that the protected attribute could be perfectly predicted from the partially repaired data set for values of $\lambda$ arbitrarily close to $1$. Thus, this upper bound provides some argument against the use of this kind of repair since the reparation should favor small distance between these two distributions to ensure a certain desired level of fairness.\vskip .1in

Hence, rather than using a displacement along the Wasserstein geodesic between the distributions, we propose the following approach called \textit{Random Repair}, that enables a better control of their Total Variation distance. \\
\indent Let $Z$ be a target variable with general distribution $\mu$ and let $B$ be a Bernoulli  variable with parameter $\lambda$, $B \sim \mathcal{B}(\lambda),$ independent of $(X,S,Y)$. Note that $R_0(Z)$ and $R_1(Z)$ follow the original conditional distributions $\mu_0$ and $\mu_1$.
Let us consider the following repair procedure which consists in randomly changing the original distribution of the variables $X$  by either selecting the target distribution $\mu$ or the original conditional distributions. The choice between both is governed by the Bernoulli parameter $\lambda$. Define  for $s \in  \{0,1\}$, the repaired distributions 
\begin{equation}\label{eq:partiallyrepairdistribution}
\tilde{\mu}_{s,\lambda}=\mathcal{L}(B Z+ (1-B)R_s(Z))=\mathcal{L}(BT_s(X) + (1-B)X \mid S=s).
\end{equation}
Note that, similarly as in the Geometric Repair, for $\lambda=0$ $\tilde{\mu}_{s,0}=\mathcal{L}(X \mid S=s)$ and for  $\lambda=1$ $\tilde{\mu}_{s,1}=\mathcal{L}(Z)=\mu$. Unlike the previous procedure, in this setting, the parameter $\lambda$ does play a role in controlling the distance between the repaired distributions
\begin{align*}
d_{TV}(\tilde{\mu}_{0,\lambda}, \tilde{\mu}_{1,\lambda})  & \leq \p(BZ + (1-B)R_0(Z) \neq BZ + (1-B)R_1(Z))\\
& = 1 - \p(BZ + (1-B)R_0(Z) = BZ + (1-B)R_1(Z)) \\ & 
\leq 1- \p(B=1) \\
& = 1 - \lambda.
\end{align*}
Hence, this bound suggests that $\lambda$ should be close to 1 to ensure non predictability of the protected attribute. \vskip .1in

Finally, observe that the misclassification error using the Randomly Repaired data is a mixture of the two errors with the totally repaired variable $T_S(X)$ and the original $X$. Thus the use of the Wasserstein barycenter $Z\sim \mu_B$ is still justified.
\begin{align*}
&R(g, \tilde{X}_{\lambda}) =\p(g(\tilde{X}_{\lambda})\neq Y)\\
& = (1-\lambda)\p(g(BT_S(X)+(1-B)X)\neq Y \mid B=0) + \lambda\p(g(BT_S(X)+(1-B)X)\neq Y \mid B=1)\\
& =(1-\lambda)\p(g(X)\neq Y \mid B=0) + \lambda\p(g(T_S(X))\neq Y \mid B=1)\\
& =(1-\lambda)\p(g(X)\neq Y) + \lambda\p(g(T_S(X))\neq Y).
\end{align*}

Therefore, in the following we promote the use of Random Repair to enhance Disparate Impact while not hampering too much the efficiency of the classification. This will be studied in the following section.

\section{Numerical Analysis of Fair Correction of a database} \label{s:num}

As the distributions at hand are empirical, the existence of an optimal transport map is not guaranteed and the repair procedure in section~\ref{s:repair} that blurs the protected variable in the original data must be adapted. In this section, we propose a new algorithm to carry this out, which in practice, achieves total fairness in contrast with the existing in the literature.

\subsection{Computational aspects}

Let $\{\left(X_i, S_i, Y_i\right), i=1, \ldots, N\}$ be an observed sample of $(X,S,Y)$, and denote by $n_0$ and $n_1$ the number of instances in each protected class. For ease of exposition and without loss of generality, suppose that the observations are ordered by the value of $S$, so we can write\\
\begin{equation*}
\begin{array}{ll}
x_{0,i}:=X_i, & \text{if} \; \; s_i=0, \ i=1,\ldots,n_0,\\
x_{1,j-n_0}:=X_j, & \text{if} \; \; s_j=1, \ j=n_0 +1,\ldots,N=n_0+n_1.
\end{array}
\end{equation*}
Generally, the sizes $n_0$ and $n_1$ of the samples $\mathcal{X}_0=\{x_{0,1}, \ldots, x_{0,n_0}\}$ and $\mathcal{X}_1=\{x_{1,1}, \ldots, x_{1,n_1}\}$ are different and Monge maps may not even exist between an empirical measure to another. This happens when their weight vectors are not compatible, which is always the case when the target measure has more points than the source measure. Hence, the solution to the optimal transport problem does not correspond to finding an optimal transport map, but an optimal transport distribution. The cuadratic cost function becomes discrete as it can be written as a matrix $C=(c_{ij})$, with $c_{ij}=\left\|x_{0,i}-x_{1,j}\right\|^2, 1\leq i \leq n_0, \ 1\leq j \leq n_1$. When $\mu_{0,n}=\sum_{i=1}^{n_0}\frac{1}{n_0}\delta_{x_{0,i}}$ and $\mu_{1,n}=\sum_{j=1}^{n_1}\frac{1}{n_1}\delta_{x_{1,j}}$, the Wasserstein distance $W_2(\mu_{0,n},\mu_{1,n})$ between them is the squared root of the optimum of a net-work flow problem known as the \textit{transportation problem}. It consists in finding a matrix $\gamma \in \mathcal{M}_{n_0 \times n_1}(\R)$ which minimizes the transportation cost between the two distributions as follows
\begin{equation}\label{eq:OTD}
\left\{	\begin{array}{cl}
\min_{{\gamma}}\displaystyle\sum_{\substack{	1 \leq i \leq n_0\\
		1\leq j \leq n_1}}c_{ij}\gamma_{ij}, &\\
\text{subject to}& \gamma_{ij} \geqslant 0,\\
& \sum_{i=1}^{n_0}\gamma_{ij}=\frac{1}{n_1}, \, \text{for all} \; j,\\
& \sum_{j=1}^{n_1}\gamma_{ij}=\frac{1}{n_0}, \, \text{for all} \; i.\\
	\end{array} \right.
\end{equation}
%\begin{equation}
%\min_{\gamma \in \Pi(\mu_{0,n_0}, \mu_{1,n_1})}\displaystyle\sum_{i,j}c_{ij}\gamma_{ij},
%\end{equation}
%where $\Pi(\mu_{0,n_0}, \mu_{1,n_1})=\left\{\gamma \in \mathcal{P}\left(\mathcal{X}_0 \times \mathcal{X}_1\right)\right.$ such that $\left. \sum_{i=1}^{n_0}\gamma_{ij}=\frac{1}{n_1}, \forall j;\sum_{j=1}^{n_1}\gamma_{ij}=\frac{1}{n_0}, \forall i \right\}$ is the set of couplings between $\mu_{0,n_0}$ and $\mu_{1,n_1}$, and
%$C=(c_{ij})_{1\leq i \leq n_0; \ 1\leq j \leq n_1}$, with $c_{ij}=\left\|x_{0,i}-x_{1,j}\right\|^2$ is the cost matrix for the Euclidean distance.

If $\hat{\gamma}$ is a solution to the linear program~\eqref{eq:OTD} then, acordingly to Remark~\ref{remark:barycenter2measures}, the distribution 
$$\mu_{B,n}=\sum_{\substack{	1 \leq i \leq n_0\\
		1\leq j \leq n_1}}\hat{\gamma}_{ij}\delta_{\{\pi_0x_{0,i}+\pi_1x_{1,j}\}}$$ is a barycenter of $\mu_{0,n}$ and $\mu_{1,n}$ with respect to weights $\pi_0$ and $\pi_1$.
See \cite{cuturi2014fast} for details on the discrete Wasserstein and Optimal Transport computation.
\subsubsection{Total repair}
In practice, the implementation of the repair scheme in section~\ref{s:repair} is based on the transport matrix $\hat{\gamma}$ from $\mathcal{X}_0$ to $\mathcal{X}_1$. As we have pointed out, in this transport scheme the major difficulty comes from the fact that the sizes of these sets are different and the transport is not a one-by-one mapping. Each point in the source set could be transported (with weights) into several points of the target, or various points in the source could be moved into the same point of the target. As a consequence, we must adapt the algorithm that produces the repaired data set, denoted by $\tilde{\mathcal{X}}$. In the following, we detail two different methods, of which the first one is similar to some existing in the literature and does not achieve total fairness in the practical framework, while the second one is a novelty and does guarantee this property for the new data $\tilde{\mathcal{X}}$.

\begin{itemize}
	\item[\textbf{(A)}] On the one hand, as depicted in Figure~\ref{fig:repair} (left), each original point in $\mathcal{X}_0, \mathcal{X}_1$ is changed by a unique point given by
	\begin{align*}
	&\tilde{x}_{0,i}=\displaystyle\frac{\displaystyle\sum_{j=1}^{n_1}\gamma_{ij}\left(\pi_0x_{0,i}+\pi_1x_{1,j}\right)}{\displaystyle\sum_{j=1}^{n_1}\gamma_{ij}}= \pi_0x_{0,i}+n_0\pi_1\sum_{j=1}^{n_1}\gamma_{ij}x_{1,j}, \ 1 \leq i \leq n_0,\\
	&\tilde{x}_{1,j}=\displaystyle\frac{\displaystyle\sum_{i=1}^{n_0}\gamma_{ij}\left(\pi_0x_{0,i}+\pi_1x_{1,j}\right)}{\displaystyle\sum_{i=1}^{n_0}\gamma_{ij}}= n_1\pi_0\displaystyle\sum_{i=1}^{n_0}\gamma_{ij}x_{0,i} + \pi_1x_{1,j}, \ 1 \leq j \leq n_1.
	\end{align*}
	Doing this, the set $\tilde{\mathcal{X}}$ will be a collection of exactly $n_0+n_1$ points. This approach generalizes to higher dimensions the idea of previous works \cite{FFMSV} and \cite{johndrow2017algorithm}, who just considere the unidimensional case, where the transport is written in terms of the distribution funtions. However, in practise it generates two sets $\tilde{\mathcal{X}}_0=\{x_{0,i}, \ 1 \leqslant i \leqslant n_0\}$ and  $\tilde{\mathcal{X}}_1=\{x_{1,j}, \ 1 \leqslant j \leqslant n_1\}$ that are not the same and do not reach \eqref{eq:fairnessrequirement}.
	%conditional distributions $\mathcal{L}\left(\tilde{X} \mid S=0\right)$ and $\mathcal{L}\left(\tilde{X} \mid S=1\right)$ that are not the same and, consequently, it does not achieve full impredictability of $S$ from $\tilde{X}$. 
	\item[\textbf{(B)}] To ensure total fairness, each point cannot be changed by a unique repaired point. Hence, each point will split its mass to be transported into several modified versions, generating an extended set $\tilde{\mathcal{X}}=\tilde{\mathcal{X}}_{0} \cup \tilde{\mathcal{X}}_{1}$, which is formed by the complete distribution $\mu_{B,n}$. More precisely, as represented in Figure~\ref{fig:repair} (right), for every $1\leq i \leq n_0, \ 1 \leq j \leq n_1,$ if $\hat{\gamma}_{ij}>0$ we define two points
	\begin{equation}\label{eq:repairFI}
	\tilde{x}_{0,i,j}=\tilde{x}_{1,j,i}:=\pi_0x_{0,i}+\pi_1x_{1,j},
	\end{equation}
	and the sets
	\begin{align*}
\tilde{\mathcal{X}}_{0}:=\bigcup\limits_{1 \leq i \leq n_0}	\{\tilde{x}_{0,i,j} \ / \ \hat{\gamma}_{ij}>0, 1 \leq j \leq n_1 \}\\
\tilde{\mathcal{X}}_{1}:=\bigcup\limits_{1 \leq j \leq n_1}	\{\tilde{x}_{1,j,i} \ / \ \hat{\gamma}_{ij}>0, 1 \leq i \leq n_0 \}.
	\end{align*}
	The rebuilt distributions have sizes equal to the number of non zero elements in $\hat{\gamma}$, and each point has weight $\hat{\gamma}_{ij}$. Unlike the previous, this approach does achieve total impredictability, as it manages to produce repaired conditional distributions equally distributed. 
\end{itemize}
%\begin{figure}[h]
%	\centering
%	\begin{subfigure}{0.43\textwidth} % width of left subfigure
%		\includegraphics[width=\textwidth]{repair1bis}
%		\caption{Not impredictability} % subcaption
%	\end{subfigure}
%	\vspace{1em} % here you can insert horizontal or vertical space
%	\begin{subfigure}{0.44\textwidth} % width of right subfigure
%		\includegraphics[width=\textwidth]{repair2}
%		\caption{Full impredictability} % subcaption
%	\end{subfigure}
%	\caption{Two repairing processes} % caption for whole figure
%	\label{fig:repair}
%\end{figure}
%\begin{figure}[h]
%	\centering
%	\begin{subfigure}{0.43\textwidth} % width of left subfigure
%		\includegraphics[width=\textwidth]{repair1plantilla}
%		\caption{Not impredictability} % subcaption
%	\end{subfigure}
%	\vspace{1em} % here you can insert horizontal or vertical space
%	\begin{subfigure}{0.44\textwidth} % width of right subfigure
%		\includegraphics[width=\textwidth]{repair2plantilla}
%		\caption{Full impredictability} % subcaption
%	\end{subfigure}
%	\caption{Two repairing processes} % caption for whole figure
%	%\label{fig:repair}
%\end{figure}
\begin{example}\label{exa:repair}
	We have simulated two samples $\mathcal{X}_0$ and $\mathcal{X}_1$ of points in $\R$ of sizes $n_0=4$ and $n_1=7$, respectively,. The optimal matrix solution to the problem~\eqref{eq:OTD} is
	\begin{equation*}
	\hat{\gamma}=\left[\begin{array}{ccccccc}
	\frac{1}{7}&\frac{1}{4}-\frac{1}{7}&0&0&0&0&0\\
	0&\frac{2}{7}-\frac{1}{4}&\frac{1}{7}&\frac{1}{14}&0&0&0\\
	0&0&0&\frac{1}{14}&\frac{1}{7}&\frac{2}{7}-\frac{1}{4}&0\\
	0&0&0&0&0&\frac{1}{4}-\frac{1}{7}&\frac{1}{7}\\
	\end{array}\right]
	\end{equation*}

If $\mathcal{X}_0$ and $\mathcal{X}_1$ are realizations of $\mathcal{L}\left(X \mid S=0\right)$ and $\mathcal{L}\left(X \mid S=1\right)$, respectively, then the left part of Figure~\ref{fig:repair} represents the blurring procedure \textbf{(A)} that produces the repaired sets $\tilde{\mathcal{X}}_0=\{\tilde{x}_{0,1}, \ldots, \tilde{x}_{0,4} \}$ (rounded green points) and $\tilde{\mathcal{X}}_1=\{\tilde{x}_{1,1}, \ldots, \tilde{x}_{1,7} \}$ (squared green points). As we can observe, the two sets are clearly different and the Statistical Parity can not be reached. Otherwise, the scheme on the right carries out the procedure \textbf{(B)}, and we note that $\tilde{\mathcal{X}}_0=\tilde{\mathcal{X}}_1$.

	\begin{figure}[h]
		\centering
		\begin{overpic}[scale=.5,unit=1mm]{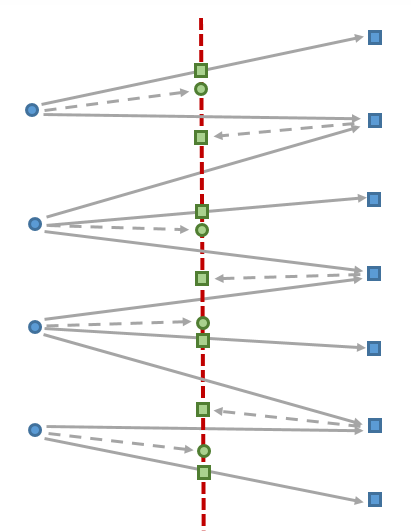}
			\put (-2,95){$S=0$}
			\put (64,98){$S=1$}
			\put (-5,78){$x_{0,1}$}
			\put (-5,57){$x_{0,2}$}
			\put (-5,37){$x_{0,3}$}
			\put (-5,17){$x_{0,4}$}
			\put (74,92){$x_{1,1}$}
			\put (74,76){$x_{1,2}$}
			\put (74,62){$x_{1,3}$}
			\put (74,48){$x_{1,4}$}
			\put (74,33){$x_{1,5}$}
			\put (74,18){$x_{1,6}$}
			\put (74,4){$x_{1,7}$}
			\put (26,85){$\tilde{x}_{0,1}$}
			\put (26,55){$\tilde{x}_{0,2}$}
			\put (26,40){$\tilde{x}_{0,3}$}
			\put (26,15){$\tilde{x}_{0,4}$}
			\put (41,90){$\tilde{x}_{1,1}$}
			\put (41,72){$\tilde{x}_{1,2}$}
			\put (41,63){$\tilde{x}_{1,3}$}
			\put (41,49){$\tilde{x}_{1,4}$}
			\put (41,32){$\tilde{x}_{1,5}$}
			\put (41,22){$\tilde{x}_{1,6}$}
			\put (41,5){$\tilde{x}_{1,7}$}
		\end{overpic}
		\hspace{45pt}
		\begin{overpic}[scale=.5,unit=1mm]{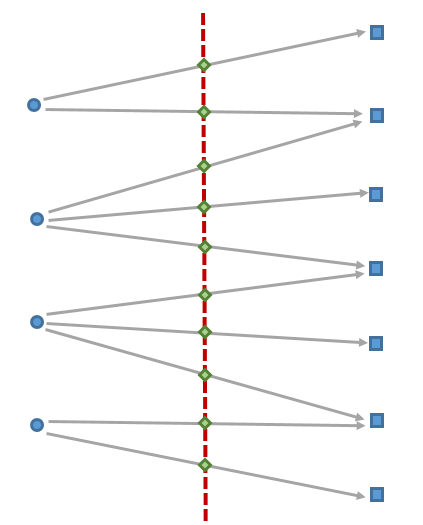}
			\put (-2,95){$S=0$}
			\put (65,99){$S=1$}
			\put (-5,78){$x_{0,1}$}
			\put (-5,57){$x_{0,2}$}
			\put (-5,37){$x_{0,3}$}
			\put (-5,17){$x_{0,4}$}
			\put (75,93){$x_{1,1}$}
			\put (75,77){$x_{1,2}$}
			\put (75,61){$x_{1,3}$}
			\put (75,48){$x_{1,4}$}
			\put (75,33){$x_{1,5}$}
			\put (75,18){$x_{1,6}$}
			\put (75,4){$x_{1,7}$}
			\put (25,90){$\tilde{x}_{0,1,1}$}
			\put (25,81){$\tilde{x}_{0,1,2}$}
			\put (25,69){$\tilde{x}_{0,2,2}$}
			\put (25,62){$\tilde{x}_{0,2,3}$}
			\put (25,50){$\tilde{x}_{0,2,4}$}
			\put (25,45){$\tilde{x}_{0,3,4}$}
			\put (25,35){$\tilde{x}_{0,3,5}$}
			\put (25,25){$\tilde{x}_{0,3,6}$}
			\put (25,15){$\tilde{x}_{0,4,6}$}
			\put (25,8){$\tilde{x}_{0,4,7}$}
			\put (42,90){$\tilde{x}_{1,1,1}$}
			\put (42,81){$\tilde{x}_{1,2,1}$}
			\put (42,69){$\tilde{x}_{1,2,2}$}
			\put (42,62){$\tilde{x}_{1,3,2}$}
			\put (42,54){$\tilde{x}_{1,4,2}$}
			\put (42,45){$\tilde{x}_{1,4,3}$}
			\put (42,37){$\tilde{x}_{1,5,3}$}
			\put (42,29){$\tilde{x}_{1,6,3}$}
			\put (42,21){$\tilde{x}_{1,6,4}$}
			\put (42,11){$\tilde{x}_{1,7,4}$}
		\end{overpic}
		\caption{Example of the two repairing processes when $n_0=4$ and $n_1=7$.}
		\label{fig:repair}
	\end{figure}
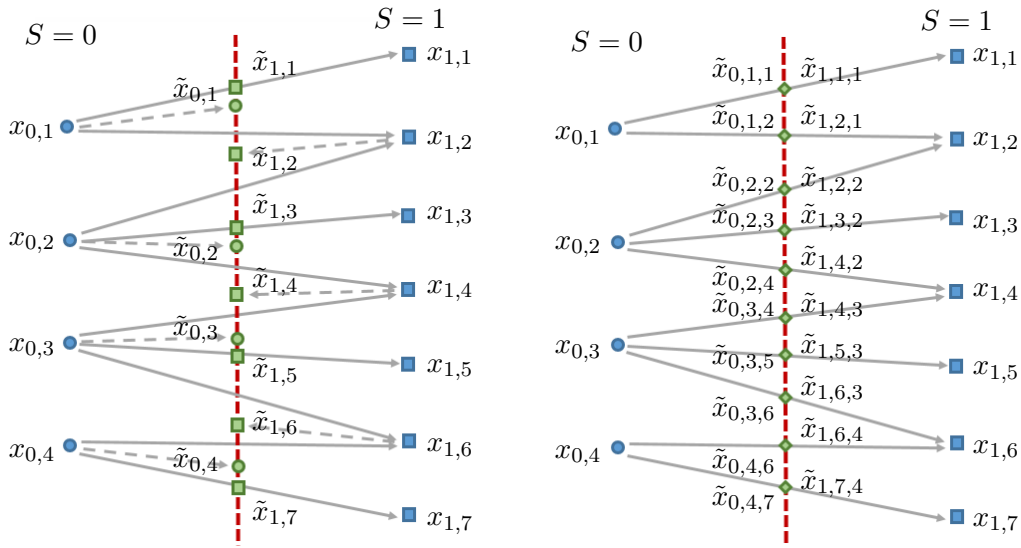
\end{example}

\begin{remark}
In the special situation when the two samples $\mathcal{X}_0$ and $\mathcal{X}_1$ have equal size $n$ and all weights are uniform, that is $\gamma_{ij}=\frac{1}{n}, \ 1 \leq i,j \leq n,$ the mass conservation constraint implies that $\gamma$ is a bijection and the Monge problem is equivalent to the optimal matching problem
\begin{equation*}
\min_{\sigma \in Perm(n)}\frac{1}{n}\sum_{i=1}^nc_{i,\sigma(i)}.
\end{equation*}
Every point in each original set will be modified by a unique point as depicted in Figure~\ref{fig:repairequalsize}. In this case, both repairing procedures \textbf{(A)} and \textbf{(B)} perfom in the same way, and total fairness is always achieved:
\begin{equation*}
\tilde{x}_{0,i}=\tilde{x}_{1,i}=\frac{1}{2}\left(x_{0,i}+x_{1,i}\right), \ 1\leq i \leq n.
\end{equation*}
\begin{figure}[H]
	\centering
	\begin{overpic}[ scale=.35,unit=1mm]{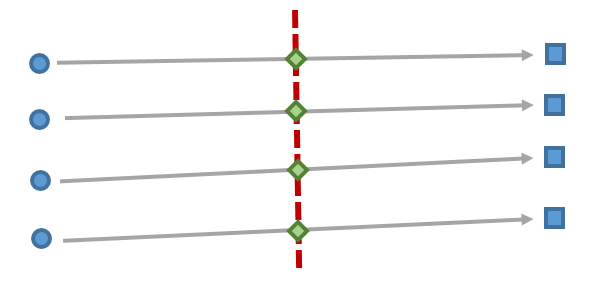}
\put (-8,38){$x_{0,1}$}
\put (-8,28){$x_{0,2}$}
\put (-8,18){$x_{0,3}$}
\put (-8,8){$x_{0,4}$}
\put (96,38){$x_{1,1}$}
\put (96,30){$x_{1,2}$}
\put (96,20){$x_{1,3}$}
\put (96,10){$x_{1,4}$}

\put (35,42){$\tilde{x}_{0,1}$}
\put (35,32){$\tilde{x}_{0,2}$}
\put (35,22){$\tilde{x}_{0,3}$}
\put (35,12){$\tilde{x}_{0,4}$}
\put (54,34){$\tilde{x}_{1,1}$}
\put (54,24){$\tilde{x}_{1,2}$}
\put (54,14){$\tilde{x}_{1,3}$}
\put (54,4){$\tilde{x}_{1,4}$}

	\end{overpic}

\caption{Repairing process when both protected groups have the same number of instances.}
\label{fig:repairequalsize}
\end{figure}
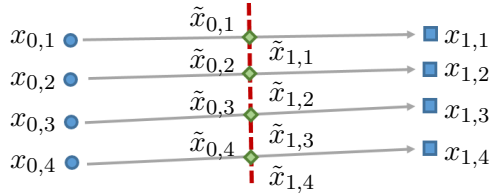
\end{remark}
\subsubsection{Random repair}
As previously noted, trying to build the set $\tilde{\mathcal{X}}$ satisfying the goal \eqref{eq:fairnessrequirement} may compromise too much the accuracy of the classification with these new data. In this sense, the Random Repair procedure proposed in section~\ref{sec:partialrepair} aims at setting a tradeoff between fairness and accuracy through the parameter $\lambda$, that models the amount of repair desired.

In this section, we detail how to compute the randomly repaired set denoted by $\tilde{\mathcal{X}}_{\lambda},$ with respect to parameter $\lambda \in \left[0,1\right]$. According to \eqref{eq:partiallyrepairdistribution}, we will randomly select either the points in the original samples $\mathcal{X}_0$ and $\mathcal{X}_1$ or their repaired sequels generated with procedure \textbf{(B)} in Figure~\ref{fig:repair} (right). More precisely, consider a sample $\{b_l\}_{l=1,\ldots, n_0+n_1} \sim B(\lambda), \ \lambda \in \left[0,1\right],$ and define
\begin{align}
&\tilde{\mathcal{X}}_{0,\lambda}:=\bigcup\limits_{1 \leq i \leq n_0}R_{0,i, \lambda}  \label{eq:repairedset0}\\
&\tilde{\mathcal{X}}_{1,\lambda}:=\bigcup\limits_{1 \leq j \leq n_1}R_{1,j, \lambda},\label{eq:repairedset1}
\end{align}
where $R_{0,i,\lambda}$ and $R_{1,j,\lambda}$ are the repaired sets of points $x_{0,i}$ and $x_{1,j}$, respectively: 
%\begin{equation}
%\tilde{x}_{0,i}= (1-b_i\pi_1)x_{0,i}+n_0b_i\pi_1\sum_{j=1}^{n_1}\gamma_{ij}x_{1,j}, \ i=1,\ldots,n_0
%\end{equation}
%\begin{equation}
%\tilde{x}_{1,j}= n_1b_j\pi_0\sum_{i=1}^{n_0}\gamma_{ij}x_{0,i} + (1-b_j\pi_0)x_{1,j}
%\end{equation}

%\begin{equation}
%\tilde{x}_{0,i,j}:= (1-b_i\pi_1)x_{0,i}+b_i\pi_1x_{1,j}, \, i=1,\ldots,n_0
%\end{equation}
%\begin{equation}
%\tilde{x}_{1,j}:= b_j\pi_0x_{0,i} + (1-b_j\pi_0)x_{1,j}, \, j=1,\ldots,n_1.
%\end{equation}

\begin{align*}
&R_{0,i, \lambda}:= \left\{ \begin{array}{ll}
\{x_{0,i}\} & \text{if} \; \; b_i=0\\
\{\tilde{x}_{0,i,j} \ / \ \hat{\gamma}_{ij}>0, 1 \leq j \leq n_1 \} & \text{if} \; \; b_i=1
\end{array}\right. 
%\ \hspace{15pt}  \text{and} \ \ \tilde{\mathcal{X}}_{0,\lambda}:=\bigcup\limits_{1 \leq i \leq n_0}R_{0,i, \lambda}
\\
&R_{1,j, \lambda}:= \left\{ \begin{array}{ll}
\{x_{1,j}\} & \text{if} \; \; b_{n_0+j}=0\\
\{\tilde{x}_{1,j,i} \ / \ \hat{\gamma}_{ij}>0, 1 \leq i \leq n_0 \} & \text{if} \; \; b_{n_0+j}=1
\end{array}\right.
% \  \text{and} \ \ \tilde{\mathcal{X}}_{1,\lambda}:=\bigcup\limits_{1 \leq j \leq n_1}R_{1,j, \lambda}
\end{align*}
with $\tilde{x}_{0,i,j}$ and $\tilde{x}_{1,j,i}$ given in Equation~\eqref{eq:repairFI}, with weights $\hat{\gamma}_{i,j}$.

\begin{example}
	Consider again the situation in Example~\ref{exa:repair}. Figure~\ref{fig:Rrepair} represents the Random Repair procedure for $\lambda=\frac{1}{2}$. We have simulated values $b_l \sim \mathcal{B}(\frac{1}{2}), \ l=1, \ldots, n_0+n_1=11,$ and the resulting sets $R_{0,i,\lambda}, 1 \leq i \leq 4,$ and $R_{1,j,\lambda}, 1\leq j \leq 7.$ Finally, from \eqref{eq:repairedset0} and \eqref{eq:repairedset1} we have the randomly repaired sets
	\begin{align*}
	&\tilde{\mathcal{X}}_{0,\lambda}=\{x_{0,1},\tilde{x}_{0,2,2},\tilde{x}_{0,2,3},\tilde{x}_{0,2,4},x_{0,3},\tilde{x}_{0,4,6},\tilde{x}_{0,4,7}\}\\
	&\tilde{\mathcal{X}}_{1,\lambda}=\{\tilde{x}_{1,1,1},x_{1,2},\tilde{x}_{1,3,2},\tilde{x}_{1,4,2},\tilde{x}_{1,4,3},\tilde{x}_{1,5,3},x_{1,6},\tilde{x}_{1,7,4}\}
	\end{align*}
\end{example}
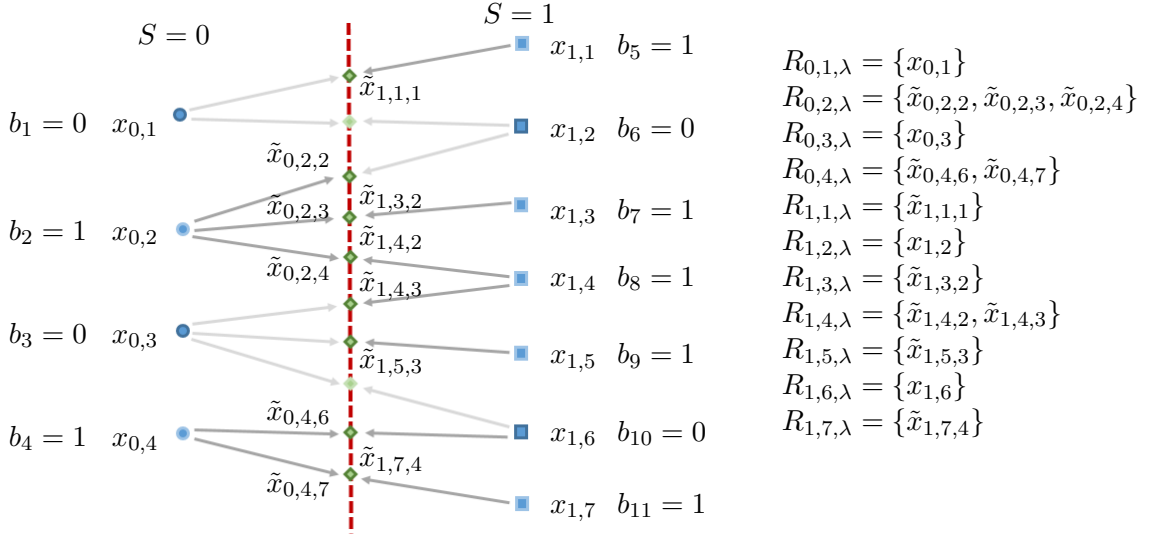
\begin{figure}[h]
\hspace{3cm}
	\begin{overpic}[scale=.5,unit=1mm]{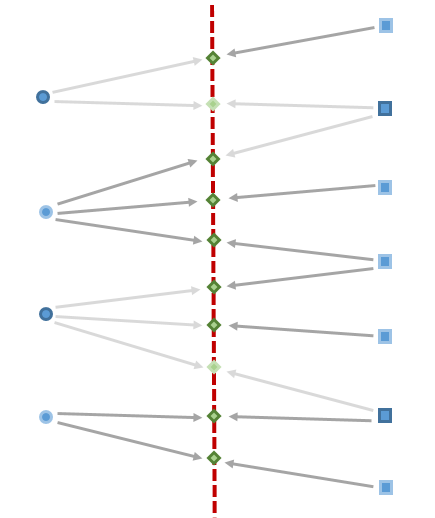}
		\put (0,95){$S=0$}
		\put (67,99){$S=1$}
		\put (-5,78){$x_{0,1}$}
		\put (-5,57){$x_{0,2}$}
		\put (-5,37){$x_{0,3}$}
		\put (-5,17){$x_{0,4}$}
		\put (-25,78){$b_1=0$}
		\put (-25,57){$b_2=1$}
		\put (-25,37){$b_3=0$}
		\put (-25,17){$b_4=1$}
		\put (80,93){$x_{1,1}$}
		\put (80,77){$x_{1,2}$}
		\put (80,61){$x_{1,3}$}
		\put (80,48){$x_{1,4}$}
		\put (80,33){$x_{1,5}$}
		\put (80,18){$x_{1,6}$}
		\put (80,4){$x_{1,7}$}
		\put (93,93){$b_5=1$}
		\put (93,77){$b_6=0$}
		\put (93,61){$b_7=1$}
		\put (93,48){$b_8=1$}
		\put (93,33){$b_9=1$}
		\put (93,18){$b_{10}=0$}
		\put (93,4){$b_{11}=1$}
		
		\put (25,72){$\tilde{x}_{0,2,2}$}
		\put (25,62){$\tilde{x}_{0,2,3}$}
		\put (25,50){$\tilde{x}_{0,2,4}$}
		\put (25,22){$\tilde{x}_{0,4,6}$}
		\put (25,8){$\tilde{x}_{0,4,7}$}
		
		\put (125,90){$R_{0,1,\lambda}=\{x_{0,1}\}$}
		\put(125,83){$R_{0,2,\lambda}=\{\tilde{x}_{0,2,2},\tilde{x}_{0,2,3},\tilde{x}_{0,2,4}\}$}
		\put (125,76){$R_{0,3,\lambda}=\{x_{0,3}\}$}
		\put(125,69){$R_{0,4,\lambda}=\{\tilde{x}_{0,4,6},\tilde{x}_{0,4,7}\}$}

		\put (43,85){$\tilde{x}_{1,1,1}$}
		\put (43,64){$\tilde{x}_{1,3,2}$}
		\put (43,56){$\tilde{x}_{1,4,2}$}
		\put (43,47){$\tilde{x}_{1,4,3}$}
		\put (43,32){$\tilde{x}_{1,5,3}$}
		\put (43,13){$\tilde{x}_{1,7,4}$}
		
		\put (125,62){$R_{1,1,\lambda}=\{\tilde{x}_{1,1,1}\}$}
		\put(125,55){$R_{1,2,\lambda}=\{x_{1,2}\}$}
		\put (125,48){$R_{1,3,\lambda}=\{\tilde{x}_{1,3,2}\}$}
		\put(125,41){$R_{1,4,\lambda}=\{\tilde{x}_{1,4,2},\tilde{x}_{1,4,3}\}$}
		\put(125,34){$R_{1,5,\lambda}=\{\tilde{x}_{1,5,3}\}$}
		\put (125,27){$R_{1,6,\lambda}=\{x_{1,6}\}$}
		\put(125,20){$R_{1,7,\lambda}=\{\tilde{x}_{1,7,4}\}$}
	\end{overpic}
	\caption{Example of the Random Repair algorithm with $\lambda=\frac{1}{2}$.}
	\label{fig:Rrepair}
\end{figure}
\subsection{Application to a real example}
To ilustrate the performance of the proposed repairing procedures in section \ref{s:repair} we consider the \textit{Adult Income} data set. It contains $29.825$ instances consisting in the values of $14$ attributes, $6$ numeric and $8$ categorical, and a categorization of each person as having an income of more or less than $50,000\$$ per year. This attribute will be the target variable in the study.  In the following, we estimate the Disparate Impact using its empirical counterpart and provide a confidence interval which was established in \cite{myfairtest}. Among the rest of the categorical attributes, we focus on the sensitive attribute \textit{Gender} (\textquotedblleft male\textquotedblright or \textquotedblleft female\textquotedblright) to be the potentially protected.
%, and Race, categorized as \textquotedblleft caucasian\textquotedblright and \textquotedblleft non caucasian\textquotedblright.
As the repairing procedures work only with the numerical attributes, to check their effectiveness we will follow the next steps:
\begin{enumerate}
	\item Split the data set into the test and the learning sample using the ratio $2.500 \, / \, 27.325$.
	\item Adjust a statistical model with the learning sample to predict the target attribute using the five numerical variables: \textit{Age, Education Level, Capital Gain, Capital Loss} and \textit{Worked hours per week}. We have trained the classifiers based on logistic regression and random forests.
	\item Predict the target for the test sample with the built model and compute the misclassification error of each rule.
	\item Apply the repair procedure in $\R^5$ to the test sample described by these numerical variables.
	%First, we will study the existence of discrimination wih respect to \textit{Sex}, and secondly with respect to \textit{Race}.
	\item Predict the target for the repaired data set with the built model and compute the misclassification error again.
\end{enumerate}

In Table~\ref{tab:models} a summary of the performance of the two classification rules considered is presented. With a confidence of $95\%$, we can say that the logit classifier has Disparate Impact at level $0.555$ and the Random Forests at $0.54$, with respect to Gender.
% We also observe that the logit classifier has Disparate Impact at level $0.646$ and the Random Forests at $0.462$, with respect to Race.
  Hence, both rules are committing discrimination with respect to this sensitive variable. Now we will see how the repairing procedures studied in section~\ref{s:repair} help in blurring the protected variable.

In Table~\ref{tab:modelsrepair} we can see that in the experiments with procedure \textbf{(A)} the estimated value for DI is not exactly $1$, as we have already anticipated. On the other hand, procedure \textbf{(B)} manages to change the data in such a way that both classification rules attain Statistical Parity. Moreover, the error in the classification done with the repaired data sets is smaller when using procedure \textbf{(B)} in the two cases. In \cite{FFMSV}, they propose a generalization to higher dimension by computing the repairing procedure for each attribute. This procedure is denoted in the table with the letter \textbf{(C)}. We see that the error is smaller than with \textbf{(A)} but still much bigger than with \textbf{(B)}. Moreover, the estimated level of Disparate Impact is not 1 but it is closer to the Statistical Parity than with procedure \textbf{(A)}.

\begin{figure}[H]
	\centering
	\begin{overpic}[scale=.6,unit=1mm]{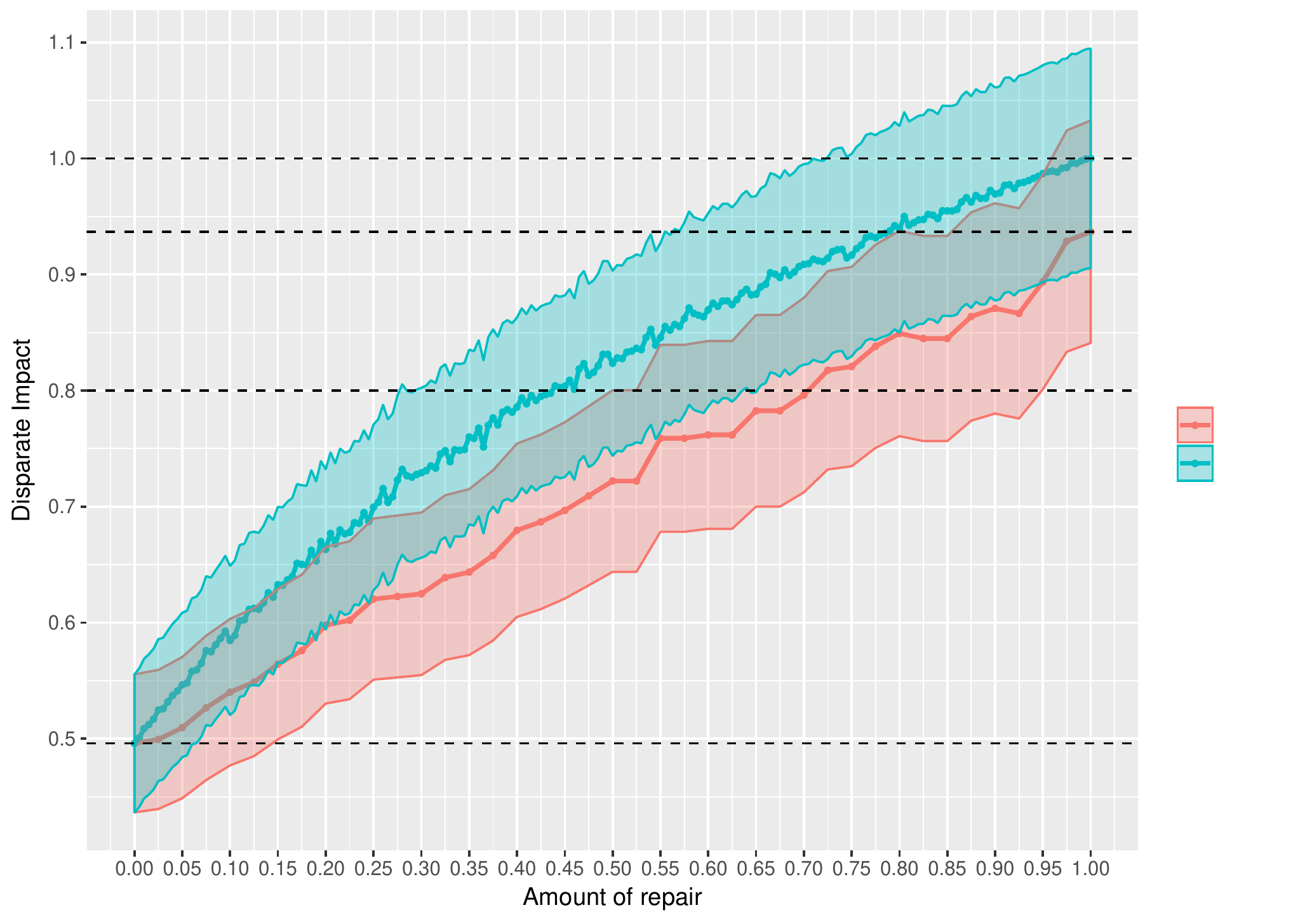}
		\put(90,42) {\scriptsize \textbf{Procedure}}
		\put(93,37) {\scriptsize Geometric}
		\put(93,34) {\scriptsize Random}
		\put(47,-1) {\scriptsize $\lambda$}
	\end{overpic}
	\caption{Confidence interval at level $95 \%$  for DI of the classifier logit with respect to Gender and the data repaired by the Geometric (red) and Random (blue) Repairs}
	\label{fig:logitSEXDI}
\end{figure}

\begin{figure}[H]
	\centering
	\begin{overpic}[scale=.6,unit=1mm]{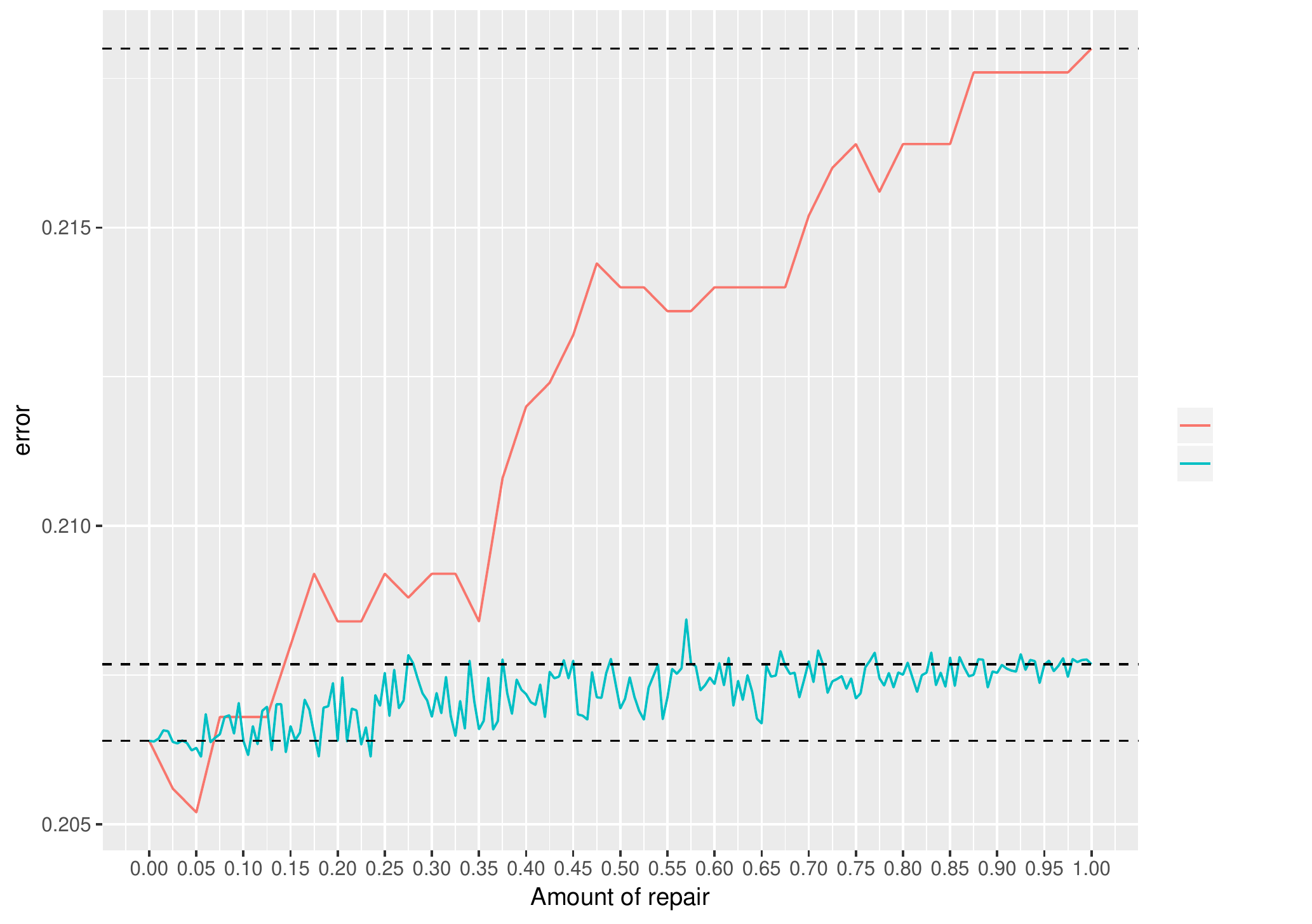}
		\put(35,72) {Classification error}
			\put(90,42) {\scriptsize \textbf{Procedure}}
		\put(93,37) {\scriptsize Geometric}
		\put(93,34) {\scriptsize Random}
		\put(47,-1) {\scriptsize $\lambda$}
	\end{overpic}
\caption{Misclassification error in the prediction with the classifier logit and the data repaired by the Geometric (red) and Random (blue) Repairs}
	\label{fig:logitSEXERROR}
\end{figure}

\begin{figure}[H]
	\centering
	\begin{overpic}[scale=.6,unit=1mm]{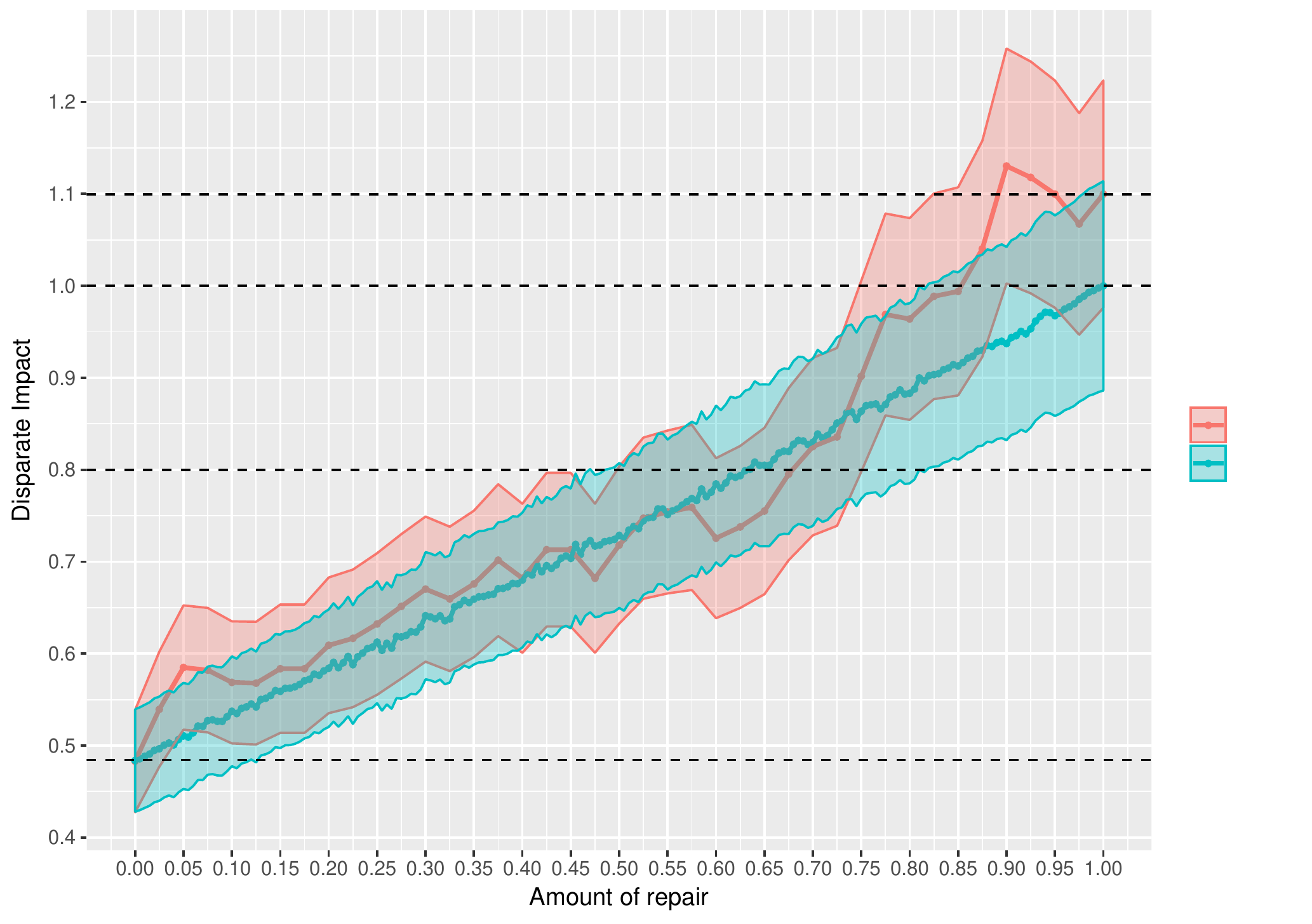}
		\put(90,42) {\scriptsize \textbf{Procedure}}
	\put(94,37) {\scriptsize Geometric}
	\put(94,34) {\scriptsize Random}
	\put(47,-1) {\scriptsize $\lambda$}
	\end{overpic}
	\caption{Confidence interval at level $95 \%$ for DI of the classifier Random Forests with respect to Gender and the data repaired by the Geometric (red) and Random (blue) Repairs}
	\label{fig:RFSEXDI}
\end{figure}

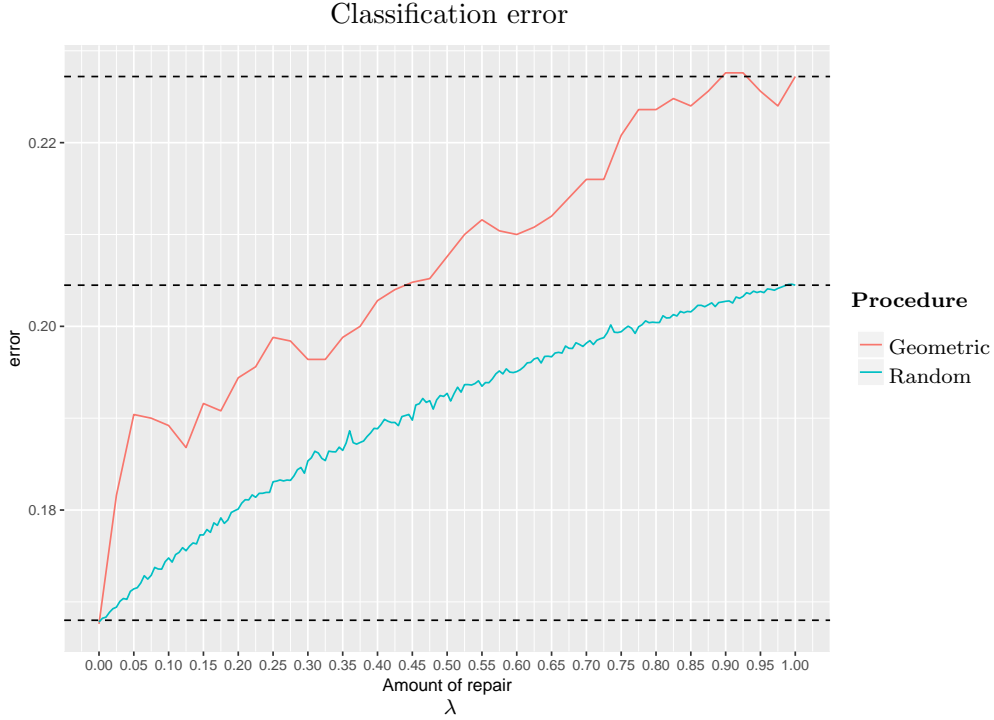
\begin{figure}[H]
	\centering
	\begin{overpic}[scale=.6,unit=1mm]{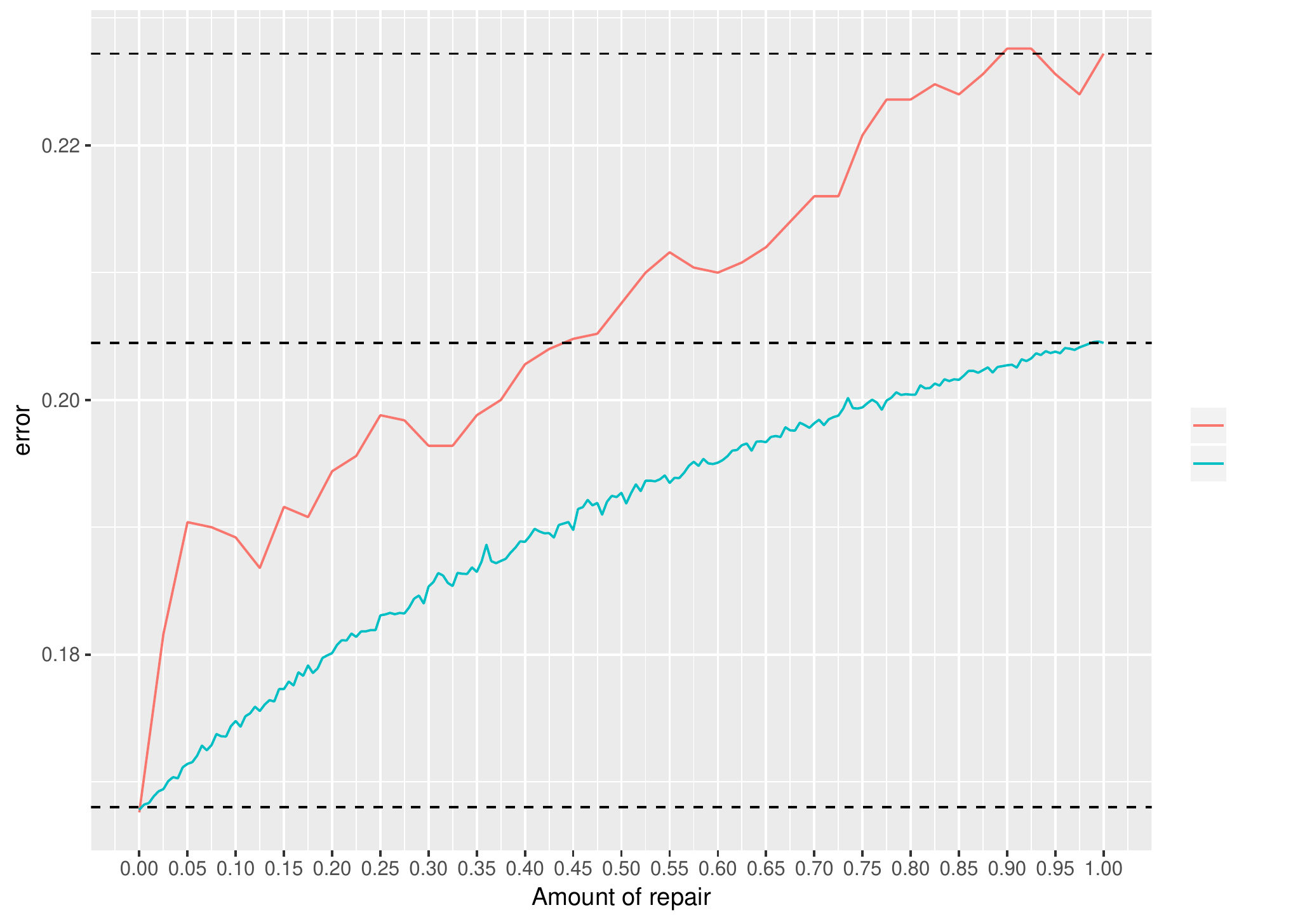}
		\put(35,72) {Classification error}
			\put(90,42) {\scriptsize \textbf{Procedure}}
		\put(94,37) {\scriptsize Geometric}
		\put(94,34) {\scriptsize Random}
		\put(47,-1) {\scriptsize $\lambda$}
	\end{overpic}
	\caption{Misclassification error in the prediction with the classifier Random Forests and the data repaired by the Geometric (red) and Random (blue) Repairs}
	\label{fig:RFSEXERROR}
\end{figure}

Finally, we present some results of the performance of the Geometric and Random Repairs. Figures~\ref{fig:logitSEXDI} and \ref{fig:RFSEXDI} represent the evolution of the estimated Disparate Impact with the amount of repair $0\leq \lambda \leq 1$. The Figures~\ref{fig:logitSEXERROR} and \ref{fig:RFSEXERROR} show the evolution with $\lambda$ of the error in the classification done from the modified data set. For the experiments concerning the Random Repair procedure (denoted RR in the figures) we have repeated it 100 times, and then we have computed the mean of the simulations. Clearly, the level of DI reached is higher with the Random Repair for the logit rule.  For the random forest procedure since the rule is not linear, the difference is not as high and Disparate Impacts have similar behaviors. Yet for larger amount of repair the gap between the two different kinds of repair increases at the advantage of the Geometric Repair.\\
Moreover, the error in the prediction from the new data modified with this procedure is smaller than with the Geometric Repair. We note that the amount of repair necessary to achieve a  confidence interval for DI at  level $0.8$ for the logit rule is
\begin{itemize}
	\item $0.3$ with the Random Repair, which entails an error of $0.2068$
	\item $0.55$ with the Geometric Repair, which entails an error of $0.2136$.
\end{itemize} 
In the case of the random forests rule, this value is $0.5$ for both but the error is
\begin{itemize}
	\item $0.1927$ with the Random Repair
	\item $0.2076$ with the Geometric Repair
\end{itemize} 
\begin{table}[H]
	\centering
	\begin{tabular}{c|ccc}
		\toprule[0.5mm]
		Statistical Model & Error & $\hat{DI}$ & CI $95\%$\\ 
		\midrule[0.5mm]
		 Logit & $0.2064$ & $0.496$ & $\left(0.437,0.555\right)$\\
		Random Forests & $0.168$ & $0.484$ & $\left(0.429,0.54\right)$\\
\bottomrule[0.5mm]
\end{tabular}
\caption{Performance and Disparate Impact with respect to the protected variable Gender.}
\label{tab:models}
\end{table}

\begin{table}[H]
	\centering
	\begin{tabular}{c|c|cccc}
		\toprule[0.5mm]
	Statistical Model & Repair & Error & Difference & $\hat{DI}$ & CI $95\%$\\ 
		\midrule[0.5mm]
		Logit & \textbf{(A)} & $0.218$ & $0.0116$ & $0.937$ & $\left(0.841,1.033\right)$\\
		Logit & \textbf{(B)} & $0.2077$ & $0.00128$ & $1$ & $\left(0.905,1.095\right)$\\
	 Logit & \textbf{(C}) & $0.2132$ & $0.0068$ & $0.94$ & $\left(0.842,1.038\right)$\\
	 	\midrule[0.3mm]
	 Random Forests & \textbf{(A)} & $0.2272$ & $0.0592$ & $1.1$ & $\left(0.976,1.223\right)$\\
		Random Forests & \textbf{(B)} & $0.2045$ & $0.0365$ & $1$ & $\left(0.886,1.114\right)$\\
		Random Forests & \textbf{(C)} & $0.2152$ & $0.0472$ & $1.091$ & $\left(0.978,1.203\right)$\\
		\bottomrule[0.5mm]
	\end{tabular}
	\caption{Repairing procedures and Disparate impact of the rules with the modified dataset}
	\label{tab:modelsrepair}
\end{table}

\section{Appendix} \label{s:append}

Proof of Theorem~\ref{theo:DIlinkBER}.
\begin{proof}
	We will show that the conditions $DI(g,X,S)\leq \tau$ and $BER(g,X,S) \leq \frac{1}{2}-\frac{a(g)}{2}(\frac{1}{\tau}-1)$ are equivalent, for all $ g \in \mathcal{G}$. Indeed, given $g \in \mathcal{G}$,
	\begin{multline*}
	BER(g,X,S) \leq \frac{1}{2}-\frac{a(g)}{2}\left(\frac{1}{\tau}-1\right)= \frac{1}{2}-\frac{(\frac{1}{\tau}-1)}{2}\p(g(X)=1 \mid S=0)\\
	\Leftrightarrow \p\left(g(X)=0 \mid S=1\right)+\p\left(g(X)=1 \mid S=0\right) \leq 1-\left(\frac{1}{\tau}-1\right)\p(g(X)=1 \mid S=0)\\
	\Leftrightarrow \left(1+\left(\frac{1}{\tau}-1\right)\right)\p\left(g(X)=1 \mid S=0\right)+\p\left(g(X)=0 \mid S=1\right) \leq 1\\
	\Leftrightarrow \frac{1}{\tau}\p\left(g(X)=1 \mid S=0\right)\leq 1- \p\left(g(X)=0 \mid S=1\right)= \p\left(g(X)=1 \mid S=1\right)\\
	\Leftrightarrow DI(g,X,S)=\frac{\p\left(g(X)=1 \mid S=0\right)}{\p\left(g(X)=1 \mid S=1\right)}\leq \tau.
	\end{multline*}
\end{proof}

Proof of Theorem~\ref{theo:minBER}
\begin{proof}
For this, we denote by $f_i, i=0,1,$ the density functions of the conditioned variables $X/S=i$, respectively, whose corresponding probability measures are both supposed to be, without loss of generality, absolute continuous with respect to a measure $\mu$. In general, the misclassification error could be written as:
\begin{multline}\label{eq:minmissER}
\p(g(X)\neq S)= \p(S=0)\p\left(g(X)=1 \mid S=0\right)+\p(S=1)\p\left(g(X)=0 \mid S=1\right)=\\
\p(S=0)\displaystyle\int_{g(X)=1}f_0(x)d\mu(x)+\p(S=1)\displaystyle\int_{g(X)=0}f_1(x)d\mu(x).
\end{multline}
Now, for $s=0,1,$ we fixe the value of $\pi_s=\p(S=s),$ and from the Bayes' Formula, we know that
\begin{equation*}
\p\left(S=s | X\right)= \frac{\pi_sf_s(X)}{\pi_0f_0(X) + \pi_1f_1(X)}.
\end{equation*}

%Suppose that we want to classify an individual for which the value of $X$ is $x$. From the Bayes' Formula, we know that:
%\begin{equation}
%\p\left(S=0 | X=x\right)= \frac{\p(S=0)\p\left(X=x | S=0\right)}{\p(S=0)\p\left(X=x | S=0\right) + \p(S=1)\p\left(X=x | S=1\right)}
%\end{equation}
%and
%\begin{equation}
%\p\left(S=1 | X=x\right)= \frac{\p(S=1)\p\left(X=x | S=1\right)}{\p(S=0)\p\left(X=x | S=0\right) + \p(S=1)\p\left(X=x | S=1\right)}.
%\end{equation}
%Now, we fixe the value of $\pi_0=\p(S=0)$ and $\pi_1=\p(S=1)$, and write these conditioned probabilities in terms of the density functions:
%\begin{equation}
%\p\left(S=0 | X\right)= \frac{\pi_0f_0(X)}{\pi_0f_0(X) + \pi_1f_1(X)}
%\end{equation}
%and
%\begin{equation}
%\p\left(S=1 | X\right)= \frac{\pi_1f_1(X)}{\pi_0f_0(X) + \pi_1f_1(X)}.
%\end{equation}
Hence,
\begin{equation*}
\left\{\p\left(S=0 | X\right) > \p\left(S=1 | X\right) \right\} = \left\{\pi_0f_0(X) > \pi_1f_1(X)\right\} , \ \mu-a.s.
\end{equation*}
%This means that we will make a smaller classification error if we classify that the individual belongs to the class $S=0$ than $S=1$. And the same reversely.
%On the other hand,
%\p(f(X)\neq S)=\displaystyle\int_{f(X)\neq S}f_0(x)d\mu(x)+\p(S=1)\displaystyle\int_{f(X)=1}f_0(x)d\mu(x).\\
%

Thus, we can deduce that the classifier that minimizes the missclassification error rate is
\begin{equation*}
g^*(x)=\left\lbrace \begin{array}{c c c}
1 & if & \pi_0f_0(x)\leq \pi_1f_1(x)\\
0 & if & \pi_0f_0(x) > \pi_1f_1(x)\\
\end{array} \right.,
\end{equation*}
and from equation (\ref{eq:minmissER}),
\begin{multline*}
\min_{g \in \mathcal{G}}\p(g(X)\neq S)= \displaystyle\int_{\left\lbrace\pi_0f_0(x)\leq \pi_1f_1(x)\right\rbrace}\pi_0f_0(x)d\mu(x)+ \displaystyle\int_{\left\lbrace\pi_0f_0(x) > \pi_1f_1(x)\right\rbrace}\pi_1f_1(x)d\mu(x).
\end{multline*}
In our particular case, $BER(g,X,S)=\p(g(X)\neq S)$ when considering $\pi_0=\pi_1=\frac{1}{2}$, so we have that
\begin{equation*}
g^*(x)=\left\lbrace \begin{array}{c c c}
1 & if & f_0(x)\leq f_1(x)\\
0 & if & f_0(x) > f_1(x)\\
\end{array} \right.
\end{equation*}
and
\begin{multline*}
\min_{g \in \mathcal{G}}BER(g,X,S)=BER(g^*,X,S)=\frac{1}{2}\left[\int_{f_0(x)\leq f_1(x)} f_0(x)d\mu(x) + \int_{f_0(x) > f_1(x)} f_1(x)d\mu(x)\right]\\
=\frac{1}{2}\int (f_0\wedge f_1)(x)d\mu(x).
\end{multline*}
This concludes the proof since  by definition
\begin{equation*}
d_{TV}(\mu_0, \mu_1)=\frac{1}{2}\int \left|f_0-f_1\right|d\mu = 1 -\int (f_0\wedge f_1)(x)d\mu(x).
\end{equation*}
\end{proof}

Proof of Lemma~\eqref{lem:bounddiffrisk}
\begin{proof}
We want to be able to control the difference $\inf_{h \in \mathcal{G}}R(h,\tilde{X}) - \inf_{g \in \mathcal{G}}R(g,X,S).$
%\begin{equation}\label{equ:differenceloss}
%\inf_{h \in \mathcal{G}}R(h,\tilde{X}) - \inf_{g \in \mathcal{G}}R(g,X,S)= R_B(\tilde{X})-R_B(X,S).
%\end{equation}
To do this, observe that
%
%\begin{multline*}
%\inf_{g \in \mathcal{G}}R(g,\tilde{X},S) \leq R(g_B, \tilde{X}, S) = R(g_B \circ T_S, X, S)=E \left[ \mathbbm{1}_{g_B\circ T_S(X)=0}(2\eta_S(X) -1)\right]+ \E \left[1-\eta_S(X) \right]
%\end{multline*}
%Thus, the difference in Equation (\ref{equ:differenceloss})
\begin{multline*}
 R_B(\tilde{X})-R_B(X,S):= \inf_{h \in \mathcal{G}}R(h,\tilde{X}) - \inf_{g \in \mathcal{G}}R(g,X,S)\\
  \leq R(g_B \circ T_S, X) - R(g_B, X, S)
  = E \left[ (2\eta_S(X) -1)(\mathbbm{1}_{g_B\circ T_S(X)=0}-\mathbbm{1}_{g_B(X,S)=0})\right]\\
 =  \E \left[(2\eta_S(X) -1)\mathbbm{1}_{g \circ T_S(X) \neq g_B (X,S)}(\mathbbm{1}_{g_B \circ T_S(X) \neq 1} - \mathbbm{1}_{g_B(X,S) \neq 1}) \right] ,
\end{multline*}
where the last equality holds because $\left( \mathbbm{1}_{g_B \circ T_S(X) \neq 1} \right) - \left(\mathbbm{1}_{g_B(X,S) \neq 1} \right) =0$ if, and only if, both classifiers have the same response $g_B \circ T_S(X) = g_B(X,S)$.

Consider $X=x$ and $S=s$,
\begin{itemize}
	\item if $g_B(x,s)=1$, $2\eta_s(x) -1 \geqslant 0$ and $\mathbbm{1}_{g_B(x,s) \neq 1}=0$. In this situation, we deduce that 
	\begin{equation*}
	\mathbbm{1}_{g_B \circ T_s(x) \neq g_B(x,s)} =1 \Leftrightarrow g_B\circ T_s(x) =0,
	\end{equation*}
	and
	\begin{equation*}
	\mathbbm{1}_{g_B \circ T_s(x) \neq 1} - \mathbbm{1}_{g_B(x,s) \neq 1}=1.
	\end{equation*}
	\item if $g_B(x,s)=0$, $2\eta_s(x) -1 <0$ and $\mathbbm{1}_{g_B(x,s) \neq 1}=1$. We deduce that 
		\begin{equation*}
	\mathbbm{1}_{g_B \circ T_s(x) \neq g_B(x,s)} =1 \Leftrightarrow g_B\circ T_s(x) =1,
	\end{equation*}
	and
	\begin{equation*}
	\mathbbm{1}_{g_B \circ T_s(x) \neq 1} - \mathbbm{1}_{g_B(x,s) \neq 1}=-1.
	\end{equation*}
\end{itemize}
In any case, the random variable $(2\eta_S(X) -1)\mathbbm{1}_{g \circ T_S(X) \neq g_B (X,S)}(\mathbbm{1}_{g_B \circ T_S(X) \neq 1} - \mathbbm{1}_{g_B(X,S) \neq 1})$ is positive and so it is its expectation
$$R(g_B\circ T_S,X)-R(g_B,X,S) = \E \left[\left| 2\eta_S(X) -1 \right| \mathbbm{1}_{g \circ T_S(X) \neq g_B (X,S)} \right]\geqslant 0.$$

Moreover, notice that $g_B \circ T_s(x)=\mathbbm{1}_{\eta_s \circ T_s(x) > \frac{1}{2} }$, for all $x$, for all $s$. Hence, $g_B\circ T_s(x) \neq g_B(x,s)$ if, and only if, either $\eta_s(x) > \frac{1}{2}$ and $\eta_s \circ T_s(x) < \frac{1}{2}$ or $\eta_s(x) < \frac{1}{2}$ and $\eta_s \circ T_s(x) > \frac{1}{2}$. In both cases, 
\begin{equation*}
\left|\eta_s(x)- \eta_s \circ T_s (x)\right| = \left|\eta_s(x) - \frac{1}{2} + \frac{1}{2} - \eta_s \circ T_s(x)\right| = \left|\eta_s(x)- \frac{1}{2} \right| +  \left| \frac{1}{2} -\eta_s \circ T_s(x)\right|,
\end{equation*}
and then it is clear that
\begin{equation*}
\left|\eta_s(x) - \frac{1}{2}\right| \leq \left|\eta_s(x)- \eta_s \circ T_s(x)\right|, \ \text{for all} \ x, \text{for all} \ s.
\end{equation*}

In conclusion, the difference between the risk using the Bayes' classifier with the original variable $X,S$ and the modified version $\tilde{X}=T_S(X)$ can be bounded as follows
\begin{equation*}
R(g_B\circ T_S,X)-R(g_B,X,S)
 \leq 2\E\left[\left|\eta_S(X)- \eta_S \circ T_S(X)\right|\right].
\end{equation*}
\end{proof}

Proof of Theorem~\ref{theo:bounddiffrisk}
\begin{proof}
	First, note that $R(h,\tilde{X})= R(h, T_S(X)) \leq R(g_B,  T_S(X))=R(g_B\circ T_S,X)$. Thus, it suffices bounding the difference between the minimal risks obtained for the best classifier with input data $(X,S)$, called $g_B$, and the risk obtained with this classification rule using the input data $\tilde{X}$
	\begin{multline*}
	R(g_B\circ T_S,X)-R(g_B,X,S)
	\leq 2\E_{(X,S)}\left[\left|\eta_S(X)- \eta_S \circ T_S(X)\right|\right]\\
	= 2\left[\p(S=0)\E_{X}\left[\left|\eta_0(X)- \eta_0 \circ T_0(X)\right| \mid S=0\right] + \p(S=1)\E_{X}\left[\left|\eta_1(X)- \eta_1 \circ T_1(X)\right| \mid S=1\right] \right]\\
	=2\displaystyle\sum_{s=0,1}\pi_s\E_X\left[\left|\eta_s(X)- \eta_s \circ T_s(X)\right| \mid S=s\right].
	\end{multline*}
	Moreover, by the Lipschitz condition and noting that $a+b \leq 2^\frac{1}{2}(a^2+b^2)^\frac{1}{2}, \ for all \ a,b \in \R$, we can write
	\begin{multline*}
	R(g_B\circ T_S,X)-R(g_B,X,S) \leq 2\displaystyle\sum_{s=0,1}\pi_sK_s\E_{X}\left[\left\|X- T_s(X)\right\| \mid S=s\right]\\
	\leq 2 \sqrt{2}K\left(\sum_{s=0,1}\pi_s^2\left(\E_{X}\left[\left\|X- T_s(X)\right\| \mid S=s\right]\right)^2\right)^\frac{1}{2},
	\end{multline*}
	where $K=\max\{K_0,K_1\}$. Finally, the Cauchy-Schwarz inequality gives
	\begin{multline*}
	R(g_B\circ T_S,X)-R(g_B,X,S) \leq 2 \sqrt{2} K \left(\sum_{s=0,1}\pi_s^2\E_{X}\left[\left\|X- T_s(X)\right\|^2 \mid S=s\right]\right)^\frac{1}{2}\\
	 = 2 \sqrt{2}K \left(\sum_{s=0,1}\pi_s^2W_2^2(\mu_s, {\mu_s}_\sharp T_s)\right)^\frac{1}{2} \leq 2 \sqrt{2}K \left(\sum_{s=0,1}\pi_sW_2^2(\mu_s, {\mu_s}_\sharp T_s)\right)^\frac{1}{2}.
	\end{multline*}
\end{proof}

\newcommand{\etalchar}[1]{$^{#1}$}

\end{document}